\documentclass[10pt,reqno,a4paper]{amsart}
\usepackage{fullpage}
\usepackage{amsmath}
\usepackage{appendix}
\usepackage{amsfonts}
\usepackage{amssymb}
\usepackage{amsxtra}
\usepackage{amsthm}
\usepackage{bm}

\usepackage{graphicx}
\usepackage{mathrsfs}
\usepackage[colorlinks=true]{hyperref}
\hypersetup{urlcolor=blue, citecolor=red, linkcolor=blue}
\usepackage[usenames,dvipsnames]{color}
\usepackage[misc]{ifsym}

\allowdisplaybreaks

\numberwithin{equation}{section}

\usepackage{cite}

\newcommand{\vertiii}[1]{{\left\vert\kern-0.25ex\left\vert\kern-0.25ex\left\vert #1 \right\vert\kern-0.25ex\right\vert\kern-0.25ex\right\vert}}

\newtheorem{thm}{Theorem}[section]

\newtheorem{lem}[thm]{Lemma}
\newtheorem{rmk}[thm]{Remark}
\newtheorem{prop}[thm]{Proposition}
{ \theoremstyle{remark} }

\begin{document}

\title[Hydrostatic limit of BE system]
{The hydrostatic limit of the Beris-Edwards system\\ in dimension two
}

\author[X. Li, M.Paicu and A. Zarnescu]{Xingyu Li, Marius Paicu and Arghir Zarnescu}
\address[Xingyu Li] { BCAM, Basque Center for Applied Mathematics, Mazarredo 14, E48009 Bilbao, Bizkaia, Spain} \email{xingyuli92@gmail.com}
\address[Marius Paicu] {Universit\'e Bordeaux, Institut de Math\'ematiques de Bordeaux, F-33405 Talence 
Cedex, France
} \email{mpaicu@math.u-bordeaux1.fr}
\address[Arghir Zarnescu] { \newline BCAM, Basque Center for Applied Mathematics, Mazarredo 14, E48009 Bilbao, Bizkaia, Spain\newline 
IKERBASQUE, Basque Foundation for Science, Maria Diaz de Haro 3, 48013, Bilbao, Bizkaia, Spain\newline
Simion Stoilow Institute of Mathematics of the Romanian Academy,
P.O. Box 1-764, RO-014700 Bucharest, Romania} \email{azarnescu@bcamath.org}
\subjclass[2010]{35Q30, 76D03.}
\keywords{Beris-Edwards system, liquid crystals, Q-tensor, hydrostatic limit.}

\begin{abstract}

We study the scaled anisotropic co-rotational Beris-Edwards system modeling the hydrodynamic motion of nematic liquid crystals in dimension two. We prove the global well-posedness with small analytic data in a thin strip domain. Moreover, we justify the limit to a system involving  the hydrostatic Navier-Stokes system with analytic data and prove the convergence.
\end{abstract}
\maketitle

\section{Introduction}
The Beris-Edwards system is a widely used model for describing the  hydrodynamic motion of nematic liquid crystals. The main characteristic feature of nematic liquid crystals, the local preferred orientation of the rod-like molecules is modeled by the so-called $Q$-tensors, that in this paper will be assumed to be specific to dimension two. The configuration space of $Q$-tensors is the set of two-by-two symmetric and traceless matrices, which is 
\[
\mathcal{S}_0^{(2)}=\{Q\in \mathbb R^{2\times 2}: Q=Q^T, \text{tr}Q=0\}
\] More details about the modeling are provided in \cite{BZ}.

\smallskip
 The Beris-Edwards system couples a   dissipative parabolic system for   $Q$-tensors-valued functions, modeling nematic liquid crystal orientation fields, with a  forced Navier-Stokes  equation for the underlying fluid velocity field $u$ of the molecules. For a general incompressible nematic liquid crystal model within the  $Q$-tensor framework,  a result about global well-posedness and decay  is provided in \cite{PZa} and also, more recently, in
\cite{SS}. The partial regularity has been recently studied by \cite{DHW} while for the related simplified Ericksen-Leslie system (which was proposed by Ericksen and Leslie in 1960's, see \cite{E1}, \cite{E2} and \cite{L}),  it was studied by Lin and Liu in \cite{LL1} and \cite{LL2}.

\bigskip
In fluid mechanics, a classical reduction of the Navier-Stokes system is obtained  assuming that  the depth of the domain and the viscosity converge to zero simultaneously,  in a related way. It is based on an approximation in geophysical fluid dynamics which assumes that the horizontal scale is large compared to the vertical scale, such that the vertical pressure gradient may be given as the product of density times the gravitational acceleration. In this case, the rescaled system is not isotropic, and we need to study the resulting  anisotropic system known as the hydrostatic Navier-Stokes system. Results and further references about the hydrostatic Navier-Stokes equation can be found in \cite{PZZ}, with  a hyperbolic version of the system available  in \cite{PZ}.
 
 \smallskip
 In our case we are interested of studying the similar problem, but in the case of the Beris-Edwards system, aiming to understand the interactions between the flow and fluid near the boundary, in the simplest possible situation, in dimension two and near a flat boundary, assming two-dimensional $Q$-tensors. Our results can be interpreted as saying that the isotropic melting condition on $Q$ at the boundary (i.e. zero Dirichlet boundary data ) is imposed, through the fluid, also  in a thin layer near the boundary, i.e. one cannot have fluid induced turbulent-like behaviour near the boundary.

\smallskip 
Thus we study the Beris-Edwards system with small analytic data in a thin strip of $\mathbb R^2  $. The equations read as follows:
\begin{equation}
\label{eqbe1}
\begin{aligned}
&\partial_t\textbf{u}+(\textbf{u}\cdot\nabla)\textbf{u}+\nabla P=\varepsilon^2\Delta\textbf{u}-\varepsilon^4\nabla\cdot(\nabla Q \odot\nabla Q+(\Delta Q)\cdot Q-Q\cdot\Delta Q)\\
\end{aligned}
\end{equation}
\begin{equation}
\label{eqbe1Q}
\nabla\cdot\textbf u=0
\end{equation}
\begin{equation}
\label{eqbe1P}
\partial_tQ+\textbf{u}\cdot\nabla Q+Q\Omega-\Omega Q=\varepsilon^2 \Delta Q-a'Q-c'Q\mbox{tr}(Q^2)
\end{equation}
where
\[
\Omega:=\frac{\nabla \textbf{u}-\nabla^T \textbf{u}}{2}
\]
The system \eqref{eqbe1}-\eqref{eqbe1P}, non-dimensionalized in a manner relevant to the study of defect patters,  is the simplified version of the one in \cite{WXZ} (see also \cite{DHW}). Here $\textbf u$ denotes the fluid velocity field, $Q$ denotes the director field, namely  a function taking values in $\mathcal{S}_0^{(2)}$ and $P$ denotes the scalar pressure function which mathematically plays the rol of the Lagrange multiplier that  guarantees the divergence free condition of the velocity field $\textbf u$. Moreover, $\Omega$ denotes the antisymmetric part of the velocity gradient tensor $\nabla\textbf u$, and $a',b',c'\in \mathbb R$ are constants, $a',c'>0$.

In \cite{WXZ}, the system contains a constant $\xi\in\mathbb R$, which is a parameter measuring the ratio between the aligning and tumbing effects that the fluid exerts on the liquid crystal molecules. In this paper, we consider the co-rotational Beris-Edwards system, namely we take $\xi=0$. Moreover, the relaxation time parameter $\Gamma$ and the fluid viscosity constant $\mu$ are both take to be one, like in \cite{DHW}.
  
In this paper, we study the system \eqref{eqbe1}-\eqref{eqbe1P}
in a thin strip in dimension two, namely in $\mathcal S^{\varepsilon}:=\{(x,y)\in \mathbb R^2, 0<y<\varepsilon\}$, and
we consider the \emph{ perturbation of shear flow }case with the scaling as follows:
\begin{equation}
\label{eqbe2}
\textbf{u}(t,x,y)=\left(u(t,x,y), v(t,x,y)\right)=\left(U(t,y/\varepsilon)+\varepsilon u^\varepsilon(t,x,y/\varepsilon), \varepsilon^2  v^\varepsilon(t,x,y/\varepsilon)\right)
\end{equation}
and
\begin{equation}
\label{eqbe3}
P(t,x,y)=\varepsilon p^\varepsilon(t,x,y/\varepsilon), \quad Q_{\alpha\beta}(t,x,y)= \begin{cases}
Q^\varepsilon_{\alpha\beta}(t,x,y/\varepsilon) \quad\mbox{if} \quad\alpha,\beta= 1\quad\mbox{or}\quad\alpha,\beta=2\\
\varepsilon Q^\varepsilon_{\alpha\beta}(t,x,y/\varepsilon)  \quad\mbox{if} \quad\alpha=1,\beta=2\quad\mbox{or}\quad\alpha=2,\beta=1
\end{cases}
\end{equation}

Furthermore, we assume the non-slip boundary condition on the fluid and isotropic boundary conditions on the director, namely:

\[\textbf u|_{y=0}=\textbf u|_{y=\varepsilon}=0, \quad Q|_{y=0}=Q|_{y=\varepsilon}=0 \]
Define $\Delta_{\varepsilon}:={\varepsilon^2}\partial_x^2+\partial_y^2$, and denote the matrix $R^\varepsilon_{ij}:=\varepsilon^4(\nabla Q \odot\nabla Q+(\Delta Q)\cdot Q-Q\cdot\Delta Q)_{ij}, 1\le i,j\le 2$. Then
\begin{equation*}
\begin{aligned}
R^{\varepsilon}_{11}&=2\varepsilon^4(\partial_xQ^{\varepsilon}_{11})^2+2\varepsilon^6(\partial_xQ^{\varepsilon}_{12})^2\\
R^{\varepsilon}_{12}&=\underbrace{2\varepsilon^3\partial_xQ^{\varepsilon}_{11}\partial_yQ^{\varepsilon}_{11}+2\varepsilon^5\partial_xQ^{\varepsilon}_{12}\partial_yQ^{\varepsilon}_{12}}_{R^\varepsilon_{12,1}}\underbrace{-2\varepsilon^3\Delta_\varepsilon Q^{\varepsilon}_{12}Q^{\varepsilon}_{11}+2\varepsilon^3Q^{\varepsilon}_{12}\Delta_\varepsilon Q^{\varepsilon}_{11}}_{R^\varepsilon_{12,2}}\\
R^{\varepsilon}_{21}&=\underbrace{2\varepsilon^3\partial_xQ^{\varepsilon}_{11}\partial_yQ^{\varepsilon}_{11}+2\varepsilon^5\partial_xQ^{\varepsilon}_{12}\partial_yQ^{\varepsilon}_{12}}_{R^\varepsilon_{21,1}}+\underbrace{2\varepsilon^3\Delta_\varepsilon Q^{\varepsilon}_{12}Q^{\varepsilon}_{11}-2\varepsilon^3Q^{\varepsilon}_{12}\Delta_\varepsilon Q^{\varepsilon}_{11}}_{R^\varepsilon_{21,2}}\\
R^{\varepsilon}_{22}&=2\varepsilon^2(\partial_yQ^{\varepsilon}_{11})^2+2\varepsilon^4(\partial_yQ^{\varepsilon}_{12})^2\\
\end{aligned}
\end{equation*}
so finally, the scaled system is considered in the strip $\mathcal S:=\{(x,y)\in\mathbb R^2, 0<y<1\}$, and the equations of $u^{\varepsilon},v^{\varepsilon}$  and $Q^{\varepsilon}_{11}, Q^{\varepsilon}_{22}$ become 
\begin{equation}
\label{equ2}
\begin{cases}
\partial_t U+\varepsilon\partial_tu^{\varepsilon}+\varepsilon(U^{\varepsilon}+ \varepsilon u^{\varepsilon})\partial_xu^{\varepsilon}+\varepsilon v^{\varepsilon}\partial_y(U^{\varepsilon}+\varepsilon  u^{\varepsilon})+\varepsilon\partial_xp^{\varepsilon}=\varepsilon^3\partial_x^2u^\varepsilon+\partial_y^2(U^{\varepsilon}+\varepsilon u^\varepsilon)-\partial_xR^\varepsilon_{11}-\partial_yR^\varepsilon_{21}\\
\varepsilon^2\partial_tv^{\varepsilon}+\varepsilon^2(U^{\varepsilon}+\varepsilon u^{\varepsilon})\partial_xv^{\varepsilon}+\varepsilon^3 v^{\varepsilon}\partial_yv^{\varepsilon}+\partial_yp^{\varepsilon}=\varepsilon^4\partial_x^2v^\varepsilon+\varepsilon^2 \partial_y^2v^\varepsilon-\partial_xR^\varepsilon_{12}-\partial_yR^\varepsilon_{22}\\
\partial_xu^{\varepsilon}+\partial_yv^{\varepsilon}=0
\end{cases}
\end{equation}
 together with the boundary condition
 \[
U|_{y=0}=U|_{y=1}=0, \quad u^\varepsilon|_{y=0}=u^\varepsilon|_{y=\varepsilon}=0, \quad v^\varepsilon|_{y=0}=v^\varepsilon|_{y=\varepsilon}=0 \]
and
 \begin{equation}
 \label{q11v2}
 \begin{cases}
 \begin{aligned}
 &\partial_tQ^{\varepsilon}_{11}+(U^{\varepsilon}+\varepsilon u^{\varepsilon})\partial_xQ^{\varepsilon}_{11}+\varepsilon v^{\varepsilon}\partial_yQ^{\varepsilon}_{11}+\partial_y(U^{\varepsilon}+\varepsilon  u^{\varepsilon})Q^{\varepsilon}_{12}-\varepsilon^3\partial_xv^{\varepsilon}Q^{\varepsilon}_{12}\\&=\varepsilon^2 \partial_x^2Q^{\varepsilon}_{11}+\partial_y^2Q^{\varepsilon}_{11}-a'Q^{\varepsilon}_{11}-2c'Q^{\varepsilon}_{11}((Q^{\varepsilon}_{11})^2+2\varepsilon^2 (Q^{\varepsilon}_{12})^2)\\
  & \varepsilon\partial_tQ_{12}^{\varepsilon}+\varepsilon (U^{\varepsilon}+\varepsilon u^{\varepsilon})\partial_x Q_{12}^{\varepsilon}+\varepsilon^2 v^{\varepsilon}\partial_yQ^{\varepsilon}_{12}-\partial_y(U^{\varepsilon}/\varepsilon +u^{\varepsilon})Q^{\varepsilon}_{11}+\varepsilon^2\partial_xv^{\varepsilon}Q^{\varepsilon}_{11} \\&=\varepsilon^3 \partial_x^2Q^{\varepsilon}_{12}+\varepsilon\partial_y^2Q^{\varepsilon}_{12}-a'\varepsilon Q^{\varepsilon}_{12}-2c'\varepsilon Q^{\varepsilon}_{12}((Q^{\varepsilon}_{11})^2+2\varepsilon^2 (Q^{\varepsilon}_{12})^2)
 \end{aligned}
 \end{cases}
 \end{equation}
 together with the boundary condition
 \[
Q_{11}^\varepsilon|_{y=0}=Q_{11}^\varepsilon|_{y=\varepsilon}=0, \quad Q_{12}^\varepsilon|_{y=0}=Q_{12}^\varepsilon|_{y=\varepsilon}=0. \]

We take  $U$ that satisfies $\partial_tU=\partial_y^2U$.
Then\footnote{ We solve $U^\varepsilon$ by separation of variables, using also the boundary conditions. Set $U^\varepsilon(t,y)=T(t)Y(y)$. Then we have $YT'=Y''T$, which is $\frac{T'}{T}=\frac{Y''}{Y}=-\gamma$ as a constant. Then we have
$T'+\gamma T=0,\quad Y''+\gamma Y=0$. We need to suppose $\gamma>0$, otherwise we will have $Y(y)\equiv 0$. By solving the ODE satistifed by $Y$, we have $Y(y)= A\cos\sqrt\gamma y+B\sin\sqrt\gamma y$.
With the boundary condition $Y(0)=Y(1)=0$, we have $\sqrt\gamma=m\pi,m\in\mathbb Z$. So $T(t)=T(0)e^{-k^2\pi^2 t}$.  We can take $k\ge 0$ without lose of  generality. }  we have
\begin{equation}
\label{shearflowu}
U(t,y)=\sum_{m\in\mathbb N}c_*(m)e^{-m^2\pi^2 t}\sin(m\pi y), \quad m\in\mathbb N.
\end{equation}

\subsection{Notations and preliminaries}\label{sec2}
We first introduce the notations that will be used in the paper as follows.
\begin{itemize}
\item[$\bullet$] We write $a\lesssim b$ to  mean that there is a uniform constant $C$, which may be different on different lines, such that $a\le Cb$.
\item[$\bullet$] We write $L^p_T(L^q_h(L^r_v))$ as the space $L^p((0,T);L^q(\mathbb R_x; L^r(\mathbb R_y)))$.
 \item[$\bullet$] We denote $(d_k)_{k\in\mathbb Z}$ (resp. $(d_k(t)))_{k\in\mathbb Z}$) to be a generic element of $l^1(\mathbb Z)$ so that $\sum_{k\in\mathbb Z} d_k=1$ (resp. $\sum_{k\in\mathbb Z} d_k(t)=1$).
  \item[$\bullet$] We denote $|D_x|$ as the Fourier multiplier with symbol $|\xi|$.
 \end{itemize}
 
\bigskip 
Next, we present some  basic aspects  of the Littlewood-Paley theory (see for \cite{BCD} for more details). For any distribution $a$ with respect to variable $x$,  we define
\begin{equation}
\label{lp1}
\Delta^h_k(a):=\mathcal F^{-1}(\phi(2^{-k}|\xi|)\hat a),\quad S^h_k(a):=\mathcal F^{-1}(\chi(2^{-k}|\xi|)\hat a)
\end{equation}
where $\phi,\chi$ are smooth functions that satisfy
\[
\mbox{Supp}\,\,\phi\subset\left\{r\in\mathbb R, \frac 34\le|r|\le\frac 83\right\}\quad \mbox{and}\quad \sum_{j\in\mathbb Z}\phi(2^{-j}r)=1\quad \mbox{for all}\quad r>0
\]

\[
\mbox{Supp}\,\,\chi\subset\left\{r\in\mathbb R, |r|\le\frac 34\right\}\quad \mbox{and}\quad \chi(r)+\sum_{j\in\mathbb Z}\phi(2^{-j}r)=1\quad \mbox{for all}\quad r>0
\]
 
 We will also need to use Bony's decomposition: 
 \begin{equation}
 \label{bony}
 fg=T^h_fg+T^h_gf+R^h(f,g)
 \end{equation}
 where
 \[
 T^f_g:=\sum_{k\in\mathbb Z}S^h_{k-1}f\Delta^h_kg, \quad R^h(f,g):=\sum_{k\in\mathbb Z}\Delta^h_kf\tilde\Delta^h_kg
 \]
 with $\tilde\Delta^h_kg:=\sum_{|k-k'|\le 1}\Delta^h_{k'}g$.  Next, we define the functional spaces as follows:
\begin{itemize}
\item[$\bullet$] For any $s\in\mathbb R$, we let
\[
||u||_{B^s}:=\sum_{k\in\mathbb Z}2^{ks}||\Delta^h_ku||_{L^2}
\]
\item[$\bullet$] For any $p\ge 1, T>0$, we define the \emph{Chemin-Lerner type} space
\[
||a||_{\tilde L^p_T(B^s)}:=\sum_{k\in\mathbb Z}2^{ks}\left(\int_0^T||\Delta^h_ka(t)||_{L^2}^pdt\right)^\frac 1p
\]
and the \emph{time weighted Chemin-Lerner type} space
\[
||a||_{\tilde L^p_{t,f}(B^s)}:=\sum_{k\in\mathbb Z}2^{ks}\left(\int_0^tf(t')||\Delta^h_ka(t')||_{L^2}^pdt'\right)^\frac 1p
\]
 \end{itemize}
 The Chemin-Lerner type space was introduced in \cite{CL} to obtain a better description of the regularizing effect of the diffusion equation. Because one cannot use a Gronwall type argument in the framework of Chemin-Lerner space, we need to use the time-weighted Chemin-Lerner norm, which was introduced in \cite{PZZ}.
 
We need an anisotropic Bernstein type lemma: \footnote{ More details can be found in \cite{CP}, \cite{P}.}
\begin{lem}
\label{bernstein}
Let $\mathcal B_h$ be a ball of $\mathbb R_h$ and $\mathcal C_h$ a ring of $\mathbb R_h$. For any $1\le p_2\le p_1\le\infty$ and $1\le q\le\infty$, we have the following results:\\

(1) If $\mbox{Supp}\,\,\hat a\subset 2^k\mathcal B_h$, then $||\partial_x^{\alpha}a||_{L^{p_1}_h(L^q_v)}\lesssim 2^{k\left(|\alpha|+\frac 1{p_2}-\frac 1{p_1}\right)}||a||_{L^{p_2}_h(L^q_v)}$.\\

(2) If $\mbox{Supp}\,\,\hat a\subset 2^k\mathcal C_h$, then $||a||_{L^{p_1}_h(L^q_v)}\lesssim 2^{-kN}||\partial^N_xa||_{L^{p_1}_h(L^q_v)}$.
\end{lem}

\bigskip
From now on, for any constant $a>0$, we define the function
\begin{equation}
\label{funpsi}
u_{\psi}(t,x,y):=\mathcal F^{-1}_{\xi\to x}(e^{\psi(t,\xi)}\hat{u}(t,\xi,y))
\end{equation}
where the phase function $\psi$ is defined by
\[
\psi(t,\xi):=(a-\lambda\eta(t))|\xi|
\]

 The function $\eta:\mathbb{R}_+\to \mathbb{R}$ and the constant $\lambda>0$ are defined in terms of $U$. Recall that $U(t,y)=\sum_{m\in\mathbb N}c_*(m)e^{-m^2\pi^2 t}\sin(m\pi y), \quad y\in [0,1], m\in\mathbb N$. For any fixed time $t$, $U$ can be seen as the Fourier series of $y$. Here we need that $U$ is uniformly bounded on $t,y$ by some constant $\epsilon>0$ small enough. 
 
We suppose that there is some $\mathsf c>0$ such that $\sum_{m>0}|c_m|<\mathsf c$ . We define the  function $\eta(t)$ as follows:
\[
\eta(0)=0,\quad \eta'(t)=\delta e^{-\mathcal Rt}+\sum_{m>0}|c_m|e^{-m^2\pi^2t}
\]
here $\mathcal R$ is half of  the constant of Poincare inequality for functions in $H^1_0(0,1)$, and $\delta<\frac c2$ needs to be small enough, where $c$   is a positive universal constant, determined by Lemma \ref{bernstein}. Moreover, $\lambda>0$ is also a universal constant which is large enough, and its definition can be found at the end of Section \ref{sec3}.
Obviously $\eta'(t)<\delta+\sum_{m>0}|c_m|<\delta+\mathsf c$, and $\eta $ needs to satisfy the bound $\eta(t)<\frac \delta{\mathcal R}+\sum_{m>0}\frac{\mathsf c}{m^2}<\frac a\lambda$.  The smallness of $\delta, \eta(t), \eta'(t)$ is necessary to prove the global well-posedness of the anisotropic system.

\subsection{Main result and strategies}
We first explain the main strategies and ideas of the proof.
 
 
 Similarly as in the Prandtl equation and in the hydrostatic Navier-Stokes equation, because of the nonlinear term $v\partial_yu$, we have one derivative loss in the $x$ variable in the process of energy estimates. So we need to work with analytic data to solve this problem. For the scaled anisotropic system, we have the following existence theorem for fixed $\varepsilon>0$.
 \begin{thm}
 \label{thm11}
 Let $a>0$.  There exists a constant $c_0>0$ small enough, such that if $U$ defined in \eqref{shearflowu} satisfies $\sum_{m>0}m|c_*(m)|< c_0$, and the initial data $(u_0^\varepsilon, v_0^\varepsilon, (Q_{11})_0^\varepsilon, (Q_{12})_0^\varepsilon)$ satisfies
 \begin{equation}
 \label{eq11}
 ||e^{a|D_x|}(\varepsilon u_0^\varepsilon,\varepsilon^2v_0^\varepsilon)||_{B^\frac 12}+||e^{a|D_x|}((Q_{11})_0^\varepsilon,\varepsilon (Q_{12})_0^\varepsilon)||_{B^\frac 12}+||e^{a|D_x|}(\varepsilon^2\partial_x,\varepsilon\partial_y)((Q_{11})_0^\varepsilon,\varepsilon (Q_{12})_0^\varepsilon)||_{B^\frac 12}\le c_0
 \end{equation}
then the system \eqref{eqbe1}-\eqref{eqbe2} has a unique global solution $(u^\varepsilon, v^\varepsilon, (Q_{11})^\varepsilon, (Q_{12})^\varepsilon)$ and there exists a constant $C>0$, such that for any $t>0$, there holds
 
\begin{equation}
\label{eq12}
\begin{aligned}
&||e^{\mathcal Rt'}(\varepsilon u^{\varepsilon},\varepsilon^2v^{\varepsilon})_\psi||_{\tilde L^\infty_t(B^\frac 12)}+||e^{\mathcal Rt'}(Q_{11}^{\varepsilon},\varepsilon Q^{\varepsilon}_{12})_\psi||_{\tilde L^\infty_t(B^\frac 12)}+||e^{\mathcal Rt'}(\varepsilon^2\partial_x,\varepsilon\partial_y)(Q_{11}^{\varepsilon},\varepsilon Q^{\varepsilon}_{12})_\psi||_{\tilde L^\infty_t(B^\frac 12)}\\
&+||e^{\mathcal Rt'}\partial_y(U^{\varepsilon}+\varepsilon u^{\varepsilon},\varepsilon^2v^{\varepsilon})_\psi||_{\tilde L^2_t(B^\frac 12)}+||e^{\mathcal Rt'}\partial_x(\varepsilon u^{\varepsilon},\varepsilon^2v^{\varepsilon})_\psi||_{\tilde L^2_t(B^\frac 12)}\\
&+||e^{\mathcal Rt'}\partial_y(Q^{\varepsilon}_{11},\varepsilon Q^{\varepsilon}_{12})_\psi||_{\tilde L^2_t(B^\frac 12)}+||e^{\mathcal Rt'}\varepsilon\partial_x(Q^{\varepsilon}_{11},\varepsilon Q^{\varepsilon}_{12})_\psi||_{\tilde L^2_t(B^\frac 12)}+\varepsilon||e^{\mathcal Rt'}\Delta_\varepsilon(Q^{\varepsilon}_{11},\varepsilon Q^{\varepsilon}_{12})_\psi||_{\tilde L^2_t(B^\frac 12)}\\
&\le C||e^{a|D_x|}(\varepsilon u_0^\varepsilon,\varepsilon^2v_0^\varepsilon)||_{B^\frac 12}+C||e^{a|D_x|}((Q_{11})_0^\varepsilon,\varepsilon (Q_{12})_0^\varepsilon)||_{B^\frac 12}+C||e^{a|D_x|}(\varepsilon^2\partial_x,\varepsilon\partial_y)((Q_{11})_0^\varepsilon,\varepsilon (Q_{12})_0^\varepsilon)||_{B^\frac 12}\\
\end{aligned}
\end{equation}
where $(u^\varepsilon_\psi,v^\varepsilon_\psi, (Q_{11})^\varepsilon_\psi,  (Q_{12})^\varepsilon_\psi)$ are defined as in  \eqref{funpsi}. 
 \end{thm}
 
 The key point in the proof of the theorem is that  we are able to control the lower frequency part $u$ by the higher frequency part $\partial_yu$ because the domain is a thin strip and we can use  on it the Poincar\'e inequality.  The  reason we choose Besov space instead of Sobolev is that we can get the optimal rate, which is just the constant of Poincar\'e inequality.

 \bigskip
 We next explain why we need to  consider our fluid part to be a perturbation of a  shear flow. The term $Q\Omega-\Omega Q$, when considered in the anistropic Besov space of the strip, will have a $1/\varepsilon$ part, which would not allow the convergence in the limit $\varepsilon\to 0$. To solve this problem, we consider $U^\varepsilon (t,y/\varepsilon)+\varepsilon u^\varepsilon{(t,x,y/\varepsilon)}$ instead of $\varepsilon u^\varepsilon(t,x,y/\varepsilon)$. For the $U^\varepsilon$ term, the $1/\varepsilon$ will vanish because $U^{\varepsilon}$ does not depend on $x$.

 Moreover, unlike in the case of the Navier-Stokes equation in \cite{PZZ}, we also need to consider the second order  frequency part $\Delta_
 \varepsilon Q$. This is because $\Delta_
 \varepsilon Q$ appears in the equation of $\textbf u$, and  $\partial_y^2 Q$ cannot be directly estimated, as we do not have the boundary condition for $\partial_yQ$. This problem can be solved by multiplying the equation for $Q$ with  $\Delta_
 \varepsilon Q$ where we have this term and  all  the other terms can be directly controlled.

\bigskip
 We next consider the hydrostatic approximation with small analytic data. Letting $\varepsilon\to 0$ the limit of \eqref{eqbe1},\eqref{eqbe1Q} formally is:
\begin{equation}
\label{dim2limit1}
\begin{cases}
\partial_tu+U\partial_xu+v\partial_yU+\partial_xp=\partial_y^2u\\
\partial_yp=0\\
\partial_xu+\partial_yv=0\\
\end{cases}
\end{equation}
with the boundary condition
 \[
U|_{y=0}=U|_{y=1}=0, \quad u|_{y=0}=u|_{y=1}=0, \quad v|_{y=0}=v|_{y=1}=0 \]
and the function $U$ satisfies $\partial_tU=\partial^2_y U$. Similarly as the scaled anisotropic system, we obtain that $U(t,y)$ also has the similar form
\eqref{shearflowu}.

Next, recall the boundary condition $v|_{y=0}=v|_{y=1}=0$. Integrating the equation $\partial_xu+\partial_yv=0$ (over $[0,1]$ with respect to the $y$ variable) we get $\partial_x\int_0^1{u(t,x,y) dy}=0$.  Because of the assumption on the space in which the solution is, we have $u(t,x,y)\to 0$ as $|x|\to\infty$. So finally we have the compatibility condition
\[
\int_0^1{u(t,x,y) dy}=0
\]
and \eqref{eqbe1P} becomes
 \begin{equation}
 \label{q11l2}
 \begin{cases}
 \begin{aligned}
 \partial_tQ_{11}+U\partial_xQ_{11}+\partial_yUQ_{12}
 & =\partial_y^2Q_{11}-aQ_{11}-2cQ^3_{11}\\
 &\partial_yUQ_{11}=0
 \end{aligned}
 \end{cases}
 \end{equation}
We need to suppose that $\partial_yU$ is not identically 0. Then we deduce from \eqref{q11l2} that $Q_{11}=Q_{12}=0$. The result about global well-posedness of the hydrostatic system is then as follows:
\begin{thm}
\label{thm12}
Let $a>0$ and let $U(t,y)$ be defined in \eqref{shearflowu} with $\partial_yU\not\equiv 0$. Suppose that the initial data $u_0$ satisfies the compatibility condition $\int_0^1{u_0dy}=0$.

(i)Assume $\sum_{m>0}\frac{|c_*(m)|}m< c_1$ for some $c_1$ small enough. Then the system \eqref{dim2limit1} has a unique global solution $u$ that satisfies for any
$s>0$ 
\begin{equation}
\label{thm12eq1}
\begin{aligned}
&||e^{\mathcal Rt'}u_\phi||_{\tilde L^\infty_t(B^s)}+||e^{\mathcal Rt'}\partial_yu_\phi||_{\tilde L^2_t(B^s)}\le ||e^{a|D_x|}u_0||_{B^s}\\
\end{aligned}
\end{equation}

(ii) Assume $\sum_{m>0}m|c_*(m)|< c_1$ for some $c_1$ small enough. Then the system \eqref{dim2limit1} has a unique global solution $u$, and for any $s>0$, there exists a constant $C>0$, such that
\begin{equation}
\label{thm12eq1}
\begin{aligned}
&||e^{\mathcal Rt'}\partial_yu_\phi||_{\tilde L^\infty_t(B^s)}+||e^{\mathcal Rt'}\partial^2_yu_\phi||_{\tilde L^2_t(B^s)}\le C(||e^{a|D_x|}u_0||_{B^s}+||e^{a|D_x|}u_0||_{B^{s+1}}+||e^{a|D_x|}\partial_yu_0||_{B^s})\\
&||e^{\mathcal Rt'}(\partial_tu)_\phi||_{\tilde L^2_t(B^s)}+||e^{\mathcal Rt'}\partial_yu_\phi||_{\tilde L^\infty_t(B^s)}\le C(||e^{a|D_x|}\partial_yu_0||_{B^s}+||e^{a|D_x|}u_0||_{B^{s+1}})
\end{aligned}
\end{equation} where $u_\phi$ is given by \eqref{hydro1}.
\end{thm}

\bigskip
\begin{rmk}
It is worth pointing out that in order to prove the global well-posedness of the hydrostatic case, we do not need the smallness of the initial data $u_0$ as \cite{PZZ}, because in our case the function $\eta$ does not include the norm $u,v$ as in the anisotropic case. However, the smallness of  $u_0$ is still required to prove the convergence.
\end{rmk}
Finally, we study the convergence of the scaled anisotropic system to  the limit hydrostatic system. For the vanishing viscosity of the analytical solutions of Navier-Stokes system in the half space, the local in time convergence  was been studied in \cite{SC}, and for the scaled anisotropic Navier-Stokes system to the hydrostatic Navier-Stokes system, the global in time convergence  was studied in \cite{PZZ}. Following the strategy of \cite{PZZ}, we have a theorem about the convergence, global in time, as follows:
 \begin{thm}
 \label{thm3}
 Let $a>0$. Suppose there exists $c_0>0$, such that $U$  satisfies $\sum_{m>0}m|c_*(m)|< c_0$ and  $(u_0^\varepsilon, v_0^\varepsilon, (Q_{11})_0^\varepsilon, (Q_{12})_0^\varepsilon)$ satisfy \eqref{eq11}. Moreover, suppose that $u_0$ satisfies $e^{a|D_x|}u_0\in B^\frac{1}{2}\cap B^\frac 52, e^{a|D_x|}\partial_yu_0\in B^\frac 32$, $||e^{a|D_x|}u_0||_{B^\frac{1}{2}}< c_0$, and the compatibility condition $\int_0^1{u_0 dy}=0$ is satisfied. The function $v_0$ is determined by $\partial_xu_0+\partial_yu_0=0$ and $v_0=0$ for $y=0$ or 1. Define $w_1:=u^\varepsilon-u, w_2:=v^\varepsilon-v$. Then there exist constants $C, M>0$, such that
 \begin{equation}
 \begin{aligned}
 \label{eqconvergencedim2}
 &||(\varepsilon w^1_{\Theta}, \varepsilon^2 w^2_\Theta)||_{\tilde L_t^\infty(B^\frac 12)}+||(Q_{11}^{\varepsilon},\varepsilon Q^{\varepsilon}_{12})_\Theta||_{\tilde L^\infty_t(B^\frac 12)}+||(\varepsilon^2\partial_x,\varepsilon\partial_y)(Q_{11}^{\varepsilon},\varepsilon Q^{\varepsilon}_{12})_\Theta||_{\tilde L^\infty_t(B^\frac 12)}\\
&+||\partial_y(V^{\varepsilon}+\varepsilon w^1,\varepsilon^2w^2)_\Theta||_{\tilde L^2_t(B^\frac 12)}+||e^{\mathcal Rt'}\partial_x(\varepsilon u^{\varepsilon},\varepsilon^2v^{\varepsilon})_\psi||_{\tilde L^2_t(B^\frac 12)}\\
&+||\partial_y(Q^{\varepsilon}_{11},\varepsilon Q^{\varepsilon}_{12})_\Theta||_{\tilde L^2_t(B^\frac 12)}+||\varepsilon\partial_x(Q^{\varepsilon}_{11},\varepsilon Q^{\varepsilon}_{12})_\Theta||_{\tilde L^2_t(B^\frac 12)}+\varepsilon||\Delta_\varepsilon(Q^{\varepsilon}_{11},\varepsilon Q^{\varepsilon}_{12})_\psi||_{\tilde L^2_t(B^\frac 12)}\\
&\le C||e^{a|D_x|}(\varepsilon(u^\varepsilon_0- u_0),\varepsilon^2(v_0^\varepsilon-v_0))||_{B^\frac 12}+C||e^{a|D_x|}((Q_{11})_0,\varepsilon (Q_{12})_0)||_{B^\frac 12}\\
&+C||e^{a|D_x|}(\varepsilon^2\partial_x,\varepsilon\partial_y)((Q_{11})_0,\varepsilon (Q_{12})_0)||_{B^\frac 12}+CM\varepsilon\\
\end{aligned}
\end{equation} where $u_\Theta$ is determined by \eqref{functiontheta} and a similar definition holds for $Q_\Theta$.
\end{thm}
To prove this theorem, we still need to control the difference between the scaled anisotropic and hydrostatic systems, and we use energy estimates with exponential weights in the Fourier variable. The weights still depend on time as before, but the loss of the analyticity in $x$ variable needs to be considered for both systems.

 
 For the hydrostatic limit case, $Q=0$ is mandatory for dimension 2 because of the relations between $u,v$ and $Q$. In dimension 3, there will be other equations about $Q_{ij}$ for $i=1,2$, and $Q$ does not to need to be 0, so that the computations will be different and more delicate. It will be studied in a future work.
 
 \subsection{Outline of the paper} In Section \ref{sec3} we prove the global wellposedness of the scaled anisotropic system. In Section \ref{Sechydro}, we prove the global wellposedness of  the hydrostatic system and the convergence from anisotropic system to hydrostatic system.

\section{Global well-posedness }\label{sec3}

We define the maximal time $T^*$ by
\begin{equation}
\label{maxt}
T^*:=\mbox{sup}\left\{t>0, \varepsilon ||u_\psi||_{B^\frac 32}+||\partial_yu_\psi||_{B^\frac 12}+ \varepsilon ||(Q_{11}, \varepsilon Q_{12})_\psi||_{B^\frac 32}+||\partial_y (Q_{11}, \varepsilon Q_{12})_\psi||_{B^\frac 12}\le \delta e^{-\mathcal R t}\right\}
\end{equation}

Note that $u_{\psi}, v_{\psi},(Q_{11})_{\psi},(Q_{12})_{\psi}$ satisfy the equations
\begin{equation}
\label{eqpsi1}
\begin{cases}
\varepsilon\partial_tu_{\psi}+\lambda\varepsilon\eta'(t)|D_x|u_{\psi}+\varepsilon([U+\varepsilon u]\partial_xu)_{\psi}+\varepsilon(v\partial_y[U+\varepsilon u])_{\psi}+\varepsilon\partial_xp_{\psi}
=\varepsilon^3\partial_x^2u_{\psi}+\varepsilon\partial_y^2u_{\psi}-\partial_x(R^\varepsilon_{11})_{\psi}-\partial_y(R^\varepsilon_{21})_{\psi}\\
\varepsilon ^2\partial_tv_{\psi}+\varepsilon^2 \lambda\eta'(t)|D_x|v_{\psi}+\varepsilon^2([U+\varepsilon u]\partial_xv)_{\psi}+\varepsilon^3(v\partial_yv)_{\psi}+\partial_yp_{\psi}=\varepsilon^4\partial_x^2v_{\psi}+\varepsilon^2\partial_y^2v_{\psi}-\partial_x(R^\varepsilon_{12})_{\psi}-\partial_y(R^\varepsilon_{22})_{\psi}\\
\partial_xu_{\psi}+\partial_yv_{\psi}=0\\
\end{cases}
\end{equation}
and
\begin{equation}
\label{eqpsi2}
\begin{cases}
\begin{aligned}
&\partial_t(Q_{11})_{\psi}+\lambda\eta'(t)|D_x|(Q_{11})_{\psi}+([U+\varepsilon u]\partial_x Q_{11}+ \varepsilon v\partial_y Q_{11})_{\psi}+( \partial_y[U+\varepsilon u] Q_{12}-\varepsilon^3\partial_xv Q_{12})_{\psi}\\
&=\varepsilon^2\partial_x^2(Q_{11})_{\psi}+\partial_y^2(Q_{11})_{\psi}-(a'Q_{11}+c'Q_{11}(2Q_{11}^2+2\varepsilon^2Q_{12}^2))_{\psi}\\
&\varepsilon\partial_t(Q_{12})_{\psi}+\varepsilon\lambda\eta'(t)|D_x|(Q_{12})_{\psi}+(\varepsilon [U+\varepsilon u]\partial_x Q_{12}+\varepsilon^2 v\partial_y Q_{12})_{\psi}-( \partial_y[U/\varepsilon+u] Q_{11}-\varepsilon^2\partial_xv Q_{11})_{\psi}\\
&=\varepsilon^3\partial_x^2(Q_{12})_{\psi}+\varepsilon\partial_y^2(Q_{12})_{\psi}-\varepsilon\left[a'Q_{12}+c'Q_{12}(2Q_{11}^2+2\varepsilon^2Q_{12}^2)\right]_{\psi}
\end{aligned}
\end{cases}
\end{equation}
We apply the dyadic operator $\Delta^h_k$. Notice that $\partial_xu_{\psi}+\partial_yv_{\psi}=0$, so
\[
(\partial_x\Delta^h_kp_\psi, \Delta^h_ku_{\psi})_{L^2}+(\partial_y\Delta^h_kp_\psi, \Delta^h_kv_{\psi})_{L^2}=-(\Delta^h_kp_\psi, \partial_x\Delta^h_ku_{\psi}+\partial_y \Delta^h_kv_{\psi})_{L^2}=0.
\]
next, observing that $\partial_xU=0$,we have from integration by parts, \[(\Delta^h_k(\partial_xuU)_\psi, \Delta^h_ku_\psi)_{L^2}=(\Delta^h_k(\partial_xuU)_\psi, \Delta^h_ku_\psi)_{L^2}+(\Delta^h_k(u\partial_xU)_\psi, \Delta^h_ku_\psi)_{L^2}=-(\Delta^h_k(uU)_\psi, \Delta^h_k\partial_xu_\psi)_{L^2}\] 
and similarly, $(\Delta^h_k(\partial_xvU)_\psi, \Delta^h_kv_\psi)_{L^2}=(\Delta^h_k(\partial_xQ_{11}U)_\psi, \Delta^h_k(Q_{11})_\psi)_{L^2}=(\Delta^h_k(\partial_xQ_{12}U)_\psi, \Delta^h_k(Q_{12})_\psi)_{L^2}=0$. Moreover,
\begin{equation}
\label{neu2}
(\Delta^h_k(\partial_yU Q_{11})_\psi, \Delta^h_k(Q_{12})_\psi)_{L^2}=(\Delta^h_k(\partial_yU Q_{12})_\psi, \Delta^h_k(Q_{11})_\psi)_{L^2}.\end{equation} 

Consider the equations \eqref{eqpsi1},\eqref{eqpsi2}, multiply respectively by $\varepsilon\Delta^h_ku_{\psi},\varepsilon^2\Delta^h_k v_{\psi},\Delta^h_k(Q_{11})_{\psi},\varepsilon\Delta^h_k(Q_{12})_{\psi}$ and take $L^2$ scalar product, to obtain that
\begin{equation}
\label{eqpsi3}
\begin{aligned}
&\frac12\frac{d}{dt}(||\Delta^h_k(\varepsilon u)_{\psi}||_{L^2}^2+||\Delta^h_k(\varepsilon^2 v_{\psi})||_{L^2}^2)+\lambda\varepsilon^2\eta'(t)(|D_x|\Delta^h_ku_{\psi},|\Delta^h_ku_{\psi})_{L^2}+\lambda\eta'(t)(|D_x|\Delta^h_k(\varepsilon^2 v_{\psi}),|\Delta^h_k(\varepsilon^2 v_{\psi}))_{L^2}\\
&+\varepsilon^4(||\partial_x\Delta^h_ku_{\psi}||_{L^2}^2+||\partial_x\Delta^h_k(\varepsilon v_{\psi})||_{L^2}^2)+(||\partial_y\Delta^h_k(\varepsilon u)_{\psi}||_{L^2}^2+||\partial_y\Delta^h_k(\varepsilon^2 v_{\psi})||_{L^2}^2)\\
&=-\varepsilon^2(\Delta^h_k(vU)_{\psi},\Delta^h_k\partial_yu_{\psi})_{L^2}-\varepsilon^3(\Delta^h_k(u\partial_xu)_{\psi},\Delta^h_ku_{\psi})_{L^2}-\varepsilon^3(\Delta^h_k(v\partial_y u)_{\psi},\Delta^h_ku_{\psi})_{L^2}\\&-\varepsilon^5(\Delta^h_k (u\partial_xv)_{\psi},\Delta^h_kv_{\psi})_{L^2}-\varepsilon^5(\Delta^h_k(v\partial_yv)_{\psi},\Delta^h_kv_{\psi})_{L^2}
+\varepsilon(\Delta^h_k(R^\varepsilon_{11})_{\psi},\Delta^h_k\partial_xu_{\psi})_{L^2}\\&+\varepsilon(\Delta^h_k(R^\varepsilon_{21})_{\psi},\Delta^h_k\partial_yu_{\psi})_{L^2}+\varepsilon^2(\Delta^h_k(R^\varepsilon_{12})_{\psi},\Delta^h_k\partial_xv_{\psi})_{L^2}+\varepsilon^2(\Delta^h_k(R^\varepsilon_{22})_{\psi},\Delta^h_k\partial_yv_{\psi})_{L^2}
\end{aligned}
\end{equation}
and
\begin{equation}
\label{eqpsi4}
\begin{aligned}
&\frac12\frac{d}{dt}(||\Delta^h_k (Q_{11})_{\psi}||_{L^2}^2+||\Delta^h_k(\varepsilon (Q_{12})_{\psi})||_{L^2}^2)+a'(||\Delta^h_k(Q_{11})_{\psi}||_{L^2}^2+||\Delta^h_k(\varepsilon Q_{12})_{\psi}||_{L^2}^2)\\
&+\lambda\eta'(t)(|D_x|\Delta^h_k(Q_{11})_{\psi},|\Delta^h_k(Q_{11})_{\psi})_{L^2}+\lambda\eta'(t)(|D_x|\Delta^h_k(\varepsilon (Q_{12})_{\psi}),|\Delta^h_k(\varepsilon (Q_{12})_{\psi}))_{L^2}\\
&+\varepsilon^2(||\partial_x\Delta^h_k(Q_{11})_{\psi}||_{L^2}^2+||\partial_x\Delta^h_k(\varepsilon (Q_{12})_{\psi})||_{L^2}^2)+(||\partial_y\Delta^h_k(Q_{11})_{\psi}||_{L^2}^2+||\partial_y\Delta^h_k(\varepsilon (Q_{12})_{\psi})||_{L^2}^2)\\
&=-\varepsilon(\Delta^h_k(u\partial_x (Q_{11}))_{\psi},\Delta^h_k (Q_{11})_{\psi})_{L^2}-\varepsilon (\Delta^h_k(v\partial_y (Q_{11}))_{\psi},\Delta^h_k (Q_{11})_{\psi})_{L^2}\\
&-\varepsilon^3(\Delta^h_k(u\partial_x (Q_{12}))_{\psi},\Delta^h_k (Q_{12})_{\psi})_{L^2}-\varepsilon^3(\Delta^h_k(v\partial_y (Q_{12}))_{\psi},\Delta^h_k (Q_{12})_{\psi})_{L^2}\\
&-(\Delta^h_k(\partial_yu Q_{11}-\varepsilon^2\partial_xv Q_{11})_{\psi}, \Delta^h_k(\varepsilon Q_{12})_{\psi})_{L^2}-\varepsilon(\Delta^h_k(\partial_yu Q_{12}-\varepsilon^2\partial_xv Q_{12})_{\psi}, \Delta^h_k(Q_{11})_{\psi})_{L^2}\\
&-2c'(\Delta^h_k[Q_{11}(Q_{11}^2+\varepsilon^2Q_{12}^2)], \Delta^h_k (Q_{11})_{\psi})_{L^2}-2c'\varepsilon^2(\Delta^h_k[Q_{12}(Q_{11}^2+\varepsilon^2Q_{12}^2)], \Delta^h_k (Q_{12})_{\psi})_{L^2}
\end{aligned}
\end{equation}

\begin{rmk}
Recall that the scaled model is in the strip $\mathcal S$. Because $u,v, Q_{11}, Q_{12}$ all equal to 0 on the boundary of $y$ variable, then according to Poincar\'e inequality, there exists a universal constant $\mathcal R>0$, such that \[2\mathcal R||\Delta_k^h (u,v,Q_{11},Q_{12})_{\psi}||_{L^2}^2\le ||\partial_y\Delta_k^h  (u,v,Q_{11},Q_{12})_{\psi}||_{L^2}^2.\]
Moreover, for the equation of $Q_{11}$ and $Q_{12}$, we need the assumption $a'>0$.
\end{rmk}
Recall from (2) of Lemma \ref{bernstein} that there exists $c>0$, such that
\begin{equation}
\label{constantc}
c2^{2k}||\Delta^h_k(u_\psi,\varepsilon v_\psi)||_{L^2}^2\le ||\Delta^h_k\partial_x(u_\psi,\varepsilon v_\psi)||_{L^2}^2
\end{equation}
multiply \eqref{eqpsi3} by $e^{2\mathcal Rt}$ and integrate over $[0,t]$, to get
\begin{equation}
\label{eqpsi5}
\begin{aligned}
&\frac12(||e^{\mathcal Rt'}\Delta^h_k(\varepsilon u)_{\psi}||_{L^{\infty}_t(L^2)}^2+||e^{\mathcal Rt'}\Delta^h_k(\varepsilon^2 v_{\psi})||_{L^{\infty}_t(L^2)}^2)
+\lambda 2^k\int_0^t{\eta'(t')(||e^{\mathcal Rt'}\Delta^h_k(\varepsilon u)_{\psi}||^2_{L^2}+||e^{\mathcal Rt'}\Delta^h_k\varepsilon^2 v_{\psi}||^2_{L^2})dt'}
\\
&+\frac12\int_0^t{e^{2\mathcal R t'}(||\Delta^h_k\partial_y(\varepsilon u_{\psi})||^2_{L^2}+||\Delta^h_k\partial_y(\varepsilon^2 v_{\psi})||^2_{L^2}+c\varepsilon^42^{2k}(||\Delta^h_k u_{\psi}||^2_{L^2}+||\Delta^h_k (\varepsilon v_{\psi})||^2_{L^2}))dt'}\\
&\le||e^{a|D_x|}\Delta_k^h (\varepsilon u_0)||^2_{L^2}+||e^{a|D_x|}\Delta_k^h (\varepsilon^2 v_0)||^2_{L^2}\\
&+\underbrace{\varepsilon^3\int_0^t{|(e^{\mathcal R t'}\Delta^h_k(u\partial_xu)_{\psi},e^{\mathcal R t'}\Delta^h_ku_{\psi})_{L^2}|dt'}}_{A_1}+\underbrace{\varepsilon^3\int_0^t{|(e^{\mathcal R t'}\Delta^h_k(v\partial_y u)_{\psi},e^{\mathcal R t'}\Delta^h_ku_{\psi})_{L^2}|dt'}}_{A_2}\\
&+\underbrace{\varepsilon^5\int_0^t{|(e^{\mathcal R t'}\Delta^h_k( u\partial_xv)_{\psi},e^{\mathcal R t'}\Delta^h_kv_{\psi})_{L^2}|dt'}}_{A_3}+\underbrace{\varepsilon^5\int_0^t{|(e^{\mathcal R t'}\Delta^h_k(v\partial_yv)_{\psi},e^{\mathcal R t'}\Delta^h_kv_{\psi})_{L^2}|dt'}}_{A_4}\\
&+\underbrace{\varepsilon\int_0^t{|(e^{\mathcal R t'}\Delta^h_k(R^\varepsilon_{11})_{\psi},e^{\mathcal R t'}\Delta^h_k\partial_xu_{\psi})_{L^2}|dt'}}_{A_5}+\underbrace{\varepsilon\int_0^t{|(e^{\mathcal R t'}\Delta^h_k(R^\varepsilon_{21})_{\psi},e^{\mathcal R t'}\Delta^h_k\partial_yu_{\psi})_{L^2}|dt'}}_{A_6}\\
&+\underbrace{\varepsilon^2\int_0^t{|(e^{\mathcal R t'}\Delta^h_k(R^\varepsilon_{12})_{\psi},e^{\mathcal R t'}\Delta^h_k\partial_xv_{\psi})_{L^2}|dt'}}_{A_7}+\underbrace{\varepsilon^2\int_0^t{|(e^{\mathcal R t'}\Delta^h_k(R^\varepsilon_{22})_{\psi},e^{\mathcal R t'}\Delta^h_k\partial_yv_{\psi})_{L^2}|dt'}}_{A_8}\\
&+\underbrace{\varepsilon^2\int_0^t{|(e^{\mathcal R t'}\Delta^h_k(v\partial_yU)_{\psi},e^{\mathcal R t'}\Delta^h_ku_{\psi})_{L^2}|dt'}}_{A_9}
\end{aligned}
\end{equation}
Similarly for the equation \eqref{eqpsi4}, we have
\begin{equation}
\label{eqpsi6}
\begin{aligned}
&\frac12||e^{\mathcal R't'}\Delta^h_k(Q_{11},\varepsilon Q_{12})_{\psi}||_{L^{\infty}_t(L^2)}^2+a'\int_0^t{||e^{\mathcal R't'}\Delta^h_k(Q_{11},\varepsilon (Q_{12})_{\psi})(t')||^2_{L^2}dt'}\\&+\lambda 2^k\int_0^t{\eta'(t')||e^{\mathcal R't'}\Delta^h_k(Q_{11},\varepsilon Q_{12})_{\psi}||^2_{L^2}dt'}
\\
&+\frac12\int_0^t{e^{2\mathcal R' t'}(||\Delta^h_k\partial_y(Q_{11}, \varepsilon Q_{12})_{\psi}||^2_{L^2}+c\varepsilon^22^{2k}||\Delta^h_k (Q_{11}, \varepsilon Q_{12})_{\psi}||^2_{L^2})dt'}\\
&\le||e^{a|D_x|}\Delta_k^h (Q_{11}, \varepsilon Q_{12})_0||^2_{L^2}+\underbrace{\varepsilon\int_0^t{|(e^{\mathcal R t'}\Delta^h_k(u\partial_x(Q_{11},  \varepsilon Q_{12}))_{\psi},e^{\mathcal R t'}\Delta^h_k(Q_{11},  \varepsilon Q_{12})_{\psi})_{L^2}|dt'}}_{B_1}\\&+\underbrace{\varepsilon\int_0^t{|(e^{\mathcal R t'}\Delta^h_k(v\partial_y(Q_{11},  \varepsilon Q_{12}))_{\psi},e^{\mathcal R t'}\Delta^h_k(Q_{11},  \varepsilon Q_{12})_{\psi})_{L^2}|dt'}}_{B_2}\\
&+\underbrace{\varepsilon\int_0^t{|\left(e^{\mathcal R’ t'}\Delta^h_k\left[(\varepsilon^2\partial_xv-\partial_y u)Q_{11}\right]_{\psi},e^{\mathcal R’ t'} \Delta^h_k(Q_{12})_{\psi}\right)_{L^2}|dt'}}_{B_3}\\
&+\underbrace{\varepsilon\int_0^t{|\left(e^{\mathcal R' t'}\Delta^h_k[(\partial_yu -\varepsilon^2\partial_xv) Q_{12}]_{\psi},e^{\mathcal R't'} \Delta^h_k(Q_{11})_{\psi}\right)_{L^2}|dt'}}_{B_4}\\
&+\underbrace{2c'\int_0^t{|(e^{\mathcal R' t'}\Delta^h_k[(Q_{11}^2+\varepsilon^2Q_{12}^2)(Q_{11},  \varepsilon Q_{12})]_{\psi},e^{\mathcal R' t'} \Delta^h_k (Q_{11},  \varepsilon Q_{12})_{\psi})_{L^2}|dt'}}_{B_5}\\
\end{aligned}
\end{equation}
From Lemma 3.1-3.3 of \cite{PZZ}, we have
\[
A_1+A_2+A_3+A_4\lesssim \varepsilon^3 d_k^22^{-k}||e^{\mathcal Rt'}(u,\varepsilon v)_{\psi}||^2_{\tilde L^2_{t,\eta'(t)}(B^1)},\quad B_1\lesssim \varepsilon d_k^22^{-k}||e^{\mathcal Rt'}(Q_{11}, \varepsilon Q_{12})_{\psi}||^2_{\tilde L^2_{t,\eta'(t)}(B^1)}
\]
Recalling that $v(t,x,y)=-\int_0^y{\partial_xu(t,x,y')dy'}$ and $\partial_yU=\sum_{m>0}m\pi c_me^{-m^2\pi^2t}\cos(m\pi y)\lesssim \eta'(t)$, we have
\[
A_{9}\lesssim \varepsilon^2 d_k^22^{-k}||e^{\mathcal Rt'}u_{\psi}||^2_{\tilde L^2_{t,\eta'(t)}(B^1)}
\]
Next, noticing the fact that $||Q^3||_{\mathcal B^{1/2}}\lesssim ||Q||^3_{\mathcal B^{1/2}}$,  we use Cauchy-Schwartz inequality to obtain that
\[
B_5\lesssim  d_k^22^{-k}||e^{\mathcal R't'}(Q_{11}, \varepsilon Q_{12})_{\psi}||^4_{L^2(B^{1/2})}
\]
\smallskip
For the rest of the terms, we need to calculate terms $A$ and $B$ seperately.
\subsubsection{ Estimates of $B_2, B_3, B_4$.} We first show the following lemmas.
\begin{lem}
\label{lem32}
For any $s\in (0,\frac 12]$ and any $t\le T^*$, we have
\[
\int_0^t{ (e^{\mathcal Rt'}\Delta^h_k(w\tilde w)_\psi, e^{\mathcal Rt'}\Delta^h_k(\partial_yu)_\psi)_{L^2}|dt'}\lesssim d_k^22^{-2ks}||e^{\mathcal Rt'}w_\psi||_{\tilde L^2_{t,\eta'(t)}(B^s)}||e^{\mathcal Rt'}\tilde w_\psi||_{\tilde L^2_{t,\eta'(t)}(B^s)}
\]
\end{lem}
\begin{proof}
By using Bony's decomposition, we have $w\tilde w=T^h_w\tilde w+T^h_{\tilde w}w+R^h(w,\tilde w)$.

\textbf{Step 1: Estimate of $\int_0^t{(e^{\mathcal Rt'}\Delta^h_k(T^h_w\tilde w)_{\psi}, e^{\mathcal Rt'}\Delta^h_k(\partial_yu)_\psi)_{L^2}dt'}$.} We have
\begin{equation*}
\begin{aligned}
&\int_0^t{(e^{\mathcal Rt'}\Delta^h_k(T^h_w\tilde w)_{\psi}, e^{\mathcal Rt'}\Delta^h_k(\partial_yu)_\psi)_{L^2}dt'}\\
&\lesssim \sum_{|k'-k|\le 4}\int_0^t{||e^{\mathcal Rt'}S^h_{k'-1}w_\psi(t')||_{L^\infty}||e^{\mathcal Rt'}\Delta^h_{k'}\tilde w_\psi(t')||_{L^2}||\Delta^h_k(\partial_yu)_\psi(t')||_{L^2}dt'}\\
&\lesssim 2^{-\frac{k}{2}} d_k(t)\sum_{|k'-k|\le 4}\int_0^t{||e^{\mathcal Rt'}S^h_{k'-1}w_\psi(t')||_{L^\infty}||e^{\mathcal Rt'}\Delta^h_{k'}\tilde w_\psi(t')||_{L^2}||\partial_yu_\psi(t')||_{B^\frac 12}dt'}\\
&\lesssim 2^{-\frac{k}{2}} d_k(t)\sum_{|k'-k|\le 4}\left(\int_0^t{||e^{\mathcal Rt'}S^h_{k'-1}w_\psi(t')||^2_{L^\infty}||\partial_yu_\psi(t')||_{B^\frac 12}dt'}\right)^\frac12\left(\int_0^t{||e^{\mathcal Rt'}\Delta^h_{k'}\tilde w_\psi(t')||^2_{L^2}||\partial_yu_\psi(t')||_{B^\frac 12}dt'}\right)^\frac12\\
\end{aligned}
\end{equation*}
Notice that
\begin{equation}
\label{eqlem32}
\begin{aligned}
\left(\int_0^t{||e^{\mathcal Rt'}S^h_{k'-1}w_\psi(t')||^2_{L^\infty}||\partial_yu_\psi(t')||_{B^\frac 12}dt'}\right)^\frac12
&\lesssim\sum_{l\le k'-2}2^\frac l2\left(\int_0^t{||e^{\mathcal Rt'}\Delta^h_{l}w_\psi(t')||^2_{L^2}||\partial_yu_\psi(t')||_{B^\frac 12}dt'}\right)^\frac12\\
&\lesssim\sum_{l\le k'-2}d_l2^{l(\frac 12-s)}||e^{\mathcal Rt'}w_\psi||_{\tilde L^2_{t,\eta'(t)}(B^s)}\lesssim 2^{k'(\frac 12-s)}||e^{\mathcal Rt'}w_\psi||_{\tilde L^2_{t,\eta'(t)}(B^s)}
\end{aligned}
\end{equation}
so we obtain that
\begin{equation*}
\begin{aligned}
\int_0^t{(e^{\mathcal Rt'}\Delta^h_k(T^h_w\tilde w)_{\psi}, e^{\mathcal Rt'}\Delta^h_k(\partial_yu)_\psi)_{L^2}dt'}&\lesssim\sum_{|k'-k|\le 4} d_kd_{k'}2^{-2ks}2^{(k-k')(2s-\frac 12)}||e^{\mathcal Rt'}w_\psi||_{\tilde L^2_{t,\eta'(t)}(B^s)}||e^{\mathcal Rt'}\tilde w_\psi||_{\tilde L^2_{t,\eta'(t)}(B^s)}\\&\lesssim d_k^22^{-2ks}||e^{\mathcal Rt'}w_\psi||_{\tilde L^2_{t,\eta'(t)}(B^s)}||e^{\mathcal Rt'}\tilde w_\psi||_{\tilde L^2_{t,\eta'(t)}(B^s)}
\end{aligned}
\end{equation*}
\textbf{Step 2: Estimate of $\int_0^t{(e^{\mathcal Rt'}\Delta^h_k(T^h_{\tilde w} w)_{\psi}, e^{\mathcal Rt'}\Delta^h_k(\partial_yu)_\psi)_{L^2}dt'}$.} Similarly as Step 1, just exchange $w$ and $\tilde w$.\\
\textbf{Step 3: Estimate of $\int_0^t{(e^{\mathcal Rt'}\Delta^h_k(R^h(w,\tilde w))_{\psi}, e^{\mathcal Rt'}\Delta^h_k(\partial_yu)_\psi)_{L^2}dt'}$}. We have
\begin{equation*}
\begin{aligned}
&\int_0^t{(e^{\mathcal Rt'}\Delta^h_k(R^h(w,\tilde w))_{\psi}, e^{\mathcal Rt'}\Delta^h_k(\partial_yu)_\psi)_{L^2}dt'}\\
&\lesssim 2^\frac k2\sum_{k'\ge k-3}\int_0^t{||e^{\mathcal Rt'}\tilde\Delta^h_{k'}w_\psi(t')||_{L^2_h(L^\infty_v)}||e^{\mathcal Rt'}\Delta^h_{k'}\tilde w_\psi(t')||_{L^2}||\Delta^h_k\partial_yu_\psi(t')||_{L^2}dt'}\\
&\lesssim \sum_{k'\ge k-3}\int_0^t{||e^{\mathcal Rt'}\tilde\Delta^h_{k'}w_\psi(t')||_{L^2}||e^{\mathcal Rt'}\Delta^h_{k'}\tilde w_\psi(t')||_{L^2}||\partial_yu_\psi||_{B^\frac 12}dt'}\\
&\lesssim \sum_{k'\ge k-3}\left(\int_0^t{||e^{\mathcal Rt'}\Delta^h_{k'}\tilde w_\psi(t')||^2_{L^2}||\partial_yu_\psi||_{B^\frac 12}dt'}\right)^\frac 12\left(\int_0^t{||e^{\mathcal Rt'}\tilde\Delta^h_{k'}w_\psi(t')||^2_{L^2}||\partial_yu_\psi||_{B^\frac 12}dt'}\right)^\frac 12\\
&\lesssim 2^{-2ks} ||e^{\mathcal Rt'}w_\psi||_{\tilde L^2_{t,\eta'(t)}(B^s)}||e^{\mathcal Rt'}\tilde w_\psi||_{\tilde L^2_{t,\eta'(t)}(B^s)}\left(\sum_{k'\ge k-3}d_{k'}2^{(k-k')s}\right)^2\\
&\lesssim d_k^22^{-2ks} ||e^{\mathcal Rt'}w_\psi||_{\tilde L^2_{t,\eta'(t)}(B^s)}||e^{\mathcal Rt'}\tilde w_\psi||_{\tilde L^2_{t,\eta'(t)}(B^s)}
\end{aligned}
\end{equation*}
and the lemma is proved after gathering three terms.
\end{proof}

\begin{lem}
\label{lem34}
For any $s\in (0,\frac 12]$ and any $t\le T^*$, we have
\begin{equation*}
\begin{aligned}
\int_0^t{ (e^{\mathcal Rt'}\Delta^h_k(\partial_yu w)_\psi, e^{\mathcal Rt'}\Delta^h_k\tilde w_\psi)_{L^2}|dt'}&\lesssim d_k^22^{-2ks}||e^{\mathcal Rt'}w_\psi||_{\tilde L^2_{t,\eta'(t)}(B^s)}||e^{\mathcal Rt'}\tilde w_\psi||_{\tilde L^2_{t,\eta'(t)}(B^{s})}\\
\end{aligned}
\end{equation*}
\end{lem}
\begin{proof}
By Bony's decomposition. we have $\partial_yuw=T^h_{\partial_yu}w+T^h_w\partial_yu+R^h(\partial_yu,w)$. We recall from (3.15) of \cite{PZZ} that 
\begin{equation}
\label{bernpoin1}
||\Delta^h_ku_\psi(t)||_{L^\infty}\lesssim d_j(t)||\partial_yu_\psi(t)||_{B^\frac 12}
\end{equation}
because of Bernstein lemma \ref{bernstein} and Poincar\'e inequality.\\
\textbf{Step 1: Estimate of $\int_0^t{(e^{\mathcal Rt'}\Delta^h_k(T^h_{\partial_yu} w)_{\psi},e^{\mathcal Rt'} \Delta^h_k\tilde w_\psi)_{L^2}dt'}$.} Notice that
\begin{equation*}
\begin{aligned}
&\int_0^t{(e^{\mathcal Rt'}\Delta^h_{k}(T^h_{\partial_yu} w)_{\psi}, e^{\mathcal Rt'}\Delta^h_k\tilde w_\psi)_{L^2}dt'}\\
&\lesssim\sum_{|k'-k|\le 4}\int_0^t{||S^h_{k'-1}\partial_yu_{\psi}(t')||_{L^\infty}||e^{\mathcal Rt'}\Delta^h_{k'}w_\psi(t')||_{L^2}||e^{\mathcal Rt'}\Delta^h_{k}\tilde w_\psi(t')||_{L^2}dt'}\\
&\lesssim\sum_{|k'-k|\le 4}\int_0^t{||\partial_yu_\psi(t)||_{B^\frac 12}||e^{\mathcal Rt'}\Delta^h_{k'}w_\psi(t')||_{L^2}||e^{\mathcal Rt'}\Delta^h_{k}\tilde w_\psi(t')||_{L^2}dt'}\\
&\lesssim\sum_{|k'-k|\le 4}\left(\int_0^t{||\partial_yu_\psi(t)||_{B^\frac 12}||e^{\mathcal Rt'}\Delta^h_{k'}\tilde w_\psi(t')||^2_{L^2}dt'}\right)^\frac 12\left(\int_0^t{||\partial_yu_\psi(t)||_{B^\frac 12}||e^{\mathcal Rt'}\Delta^h_{k}w_\psi(t')||_{L^2}^2dt'}\right)^\frac 12\\
&\lesssim d_k2^{-2ks}||e^{\mathcal Rt'}w_\psi||_{\tilde L^2_{t,\eta'(t)}(B^s)}||e^{\mathcal Rt'}\tilde w_\psi||_{\tilde L^2_{t,\eta'(t)}(B^{s})}\left(\sum_{|k'-k|\le 4}d_{k'}2^{(k-k')s}\right)\\
&\lesssim d_k^22^{-2ks}||e^{\mathcal Rt'}w_\psi||_{\tilde L^2_{t,\eta'(t)}(B^s)}||e^{\mathcal Rt'}\tilde w_\psi||_{\tilde L^2_{t,\eta'(t)}(B^{s})}
\end{aligned}
\end{equation*}
\textbf{Step 2: Estimate of $\int_0^t{(e^{\mathcal Rt'}\Delta^h_k(T^h_w{\partial_yu})_{\psi}, e^{\mathcal Rt'}\Delta^h_k\tilde w_\psi)_{L^2}dt'}$.} We have
\begin{equation*}
\begin{aligned}
&\int_0^t{(e^{\mathcal Rt'}\Delta^h_k(T^h_w{\partial_yu})_{\psi},e^{\mathcal Rt'} \Delta^h_k\tilde w_\psi)_{L^2}dt'}\\
&\lesssim\sum_{|k'-k|\le 4}\int_0^t{||S^h_{k'-1}w_{\psi}(t')||_{L^\infty}||e^{\mathcal Rt'}\Delta^h_{k'}\partial_yu_\psi(t')||_{L^2}||e^{\mathcal Rt'}\Delta^h_{k}\tilde w_\psi(t')||_{L^2}dt'}\\
&\lesssim\sum_{|k'-k|\le 4}2^{-\frac{k'}{2}}\int_0^t{d_{k'}(t)||S^h_{k'-1}w_{\psi}(t')||_{L^\infty}||e^{\mathcal Rt'}\partial_yu_\psi(t')||_{B^\frac 12}||e^{\mathcal Rt'}\Delta^h_{k}\tilde w_\psi(t')||_{L^2}dt'}\\
&\lesssim \sum_{|k'-k|\le 4}d_{k'}2^{-\frac{k'}{2}}\left(\int_0^t{||S^h_{k'-1}w_{\psi}(t')||^2_{L^\infty}||e^{\mathcal Rt'}\partial_yu_\psi(t')||_{B^\frac 12}dt'}\right)^\frac 12\left(\int_0^t{||\partial_yu_\psi(t')||_{B^\frac 12}||e^{\mathcal Rt'}\Delta^h_{k}\tilde w_\psi(t')||^2_{L^2}dt'}\right)^{\frac 12}
\end{aligned}
\end{equation*}
We use \eqref{eqlem32} again to obtain that
\begin{equation*}
\begin{aligned}
\int_0^t{(e^{\mathcal Rt'}\Delta^h_k(T^h_w{\partial_yu})_{\psi},e^{\mathcal Rt'} \Delta^h_k\tilde w_\psi)_{L^2}dt'}&\lesssim\sum_{|k'-k|\le 4} d_kd_{k'}2^{-2ks}2^{(k-k')(s-\frac 12)}||e^{\mathcal Rt'}w_\psi||_{\tilde L^2_{t,\eta'(t)}(B^s)}||e^{\mathcal Rt'}\tilde w_\psi||_{\tilde L^2_{t,\eta'(t)}(B^s)}\\&\lesssim d_k^22^{-2ks}||e^{\mathcal Rt'}w_\psi||_{\tilde L^2_{t,\eta'(t)}(B^s)}||e^{\mathcal Rt'}\tilde w_\psi||_{\tilde L^2_{t,\eta'(t)}(B^s)}
\end{aligned}
\end{equation*}
\textbf{Step 3: Estimate of $\int_0^t{(e^{\mathcal Rt'}\Delta^h_k(R^h(\partial_yu,w))_{\psi},e^{\mathcal Rt'} \Delta^h_k\tilde w_\psi)_{L^2}dt'}$.} We have
\begin{equation*}
\begin{aligned}
&\int_0^t{(\Delta^h_k(e^{\mathcal Rt'}R^h(\partial_yu,w))_{\psi}, e^{\mathcal Rt'}\Delta^h_k\tilde w_\psi)_{L^2}dt'}\\
&\lesssim 2^{\frac k2}\int_0^t{||\tilde \Delta^h_{k'}\partial_yu_\psi(t')||_{L^2_h(L^\infty_v)}||e^{\mathcal Rt'}\Delta^h_{k'}w_\psi(t')||_{L^2}||e^{\mathcal Rt'}\Delta^h_k\tilde w_\psi(t')||_{L^2}dt'}\\
&\lesssim 2^{\frac k2} \sum_{k'\ge k-3}2^{-\frac {k'}2}\int_0^t{||\partial_yu_\psi(t')||_{B^\frac 12}||e^{\mathcal Rt'}\Delta^h_{k'} w_\psi(t')||_{L^2}||e^{\mathcal Rt'}\Delta^h_k\tilde w_\psi(t')||_{L^2}dt'}\\
&\lesssim 2^{\frac k2}\sum_{k'\ge k-3}2^{-\frac {k'}2}\left(\int_0^t{||e^{\mathcal Rt'}\Delta^h_{k}\tilde w_\psi(t')||^2_{L^2}||\partial_yu_\psi||_{B^\frac 12}dt'}\right)^\frac 12\left(\int_0^t{||e^{\mathcal Rt'}\Delta^h_{k'}w_\psi(t')||^2_{L^2}||\partial_yu_\psi||_{B^\frac 12}dt'}\right)^\frac 12\\
&\lesssim d_k2^{-2ks} ||e^{\mathcal Rt'}w_\psi||_{\tilde L^2_{t,\eta'(t)}(B^s)}||e^{\mathcal Rt'}\tilde w_\psi||_{\tilde L^2_{t,\eta'(t)}(B^s)}\left(\sum_{k'\ge k-3}d_{k'}2^{(k-k')(s+\frac 12)}\right)\\
&\lesssim d_k^22^{-2ks} ||e^{\mathcal Rt'}w_\psi||_{\tilde L^2_{t,\eta'(t)}(B^s)}||e^{\mathcal Rt'}\tilde w_\psi||_{\tilde L^2_{t,\eta'(t)}(B^s)}
\end{aligned}
\end{equation*}
so we finish the proof after gathering all the terms.
\end{proof}
\begin{rmk}
Lemma \ref{lem32} and Lemma \ref{lem34} give the estimates on the terms $A_6$ and the $B_3$ term with $\partial_yu$ part. But since for these terms, there is no $\partial_x$ or $v$ term included, the Besov norm becomes $B^\frac 12$ instead of $B^1$, which makes the  function $\eta(t)$  a bit different from the one in \cite{PZZ}.
\end{rmk}

\begin{prop}
\label{propb4}
\begin{equation}
\label{estib4}
B_2\lesssim \frac 1{\varepsilon} d_k^22^{-k}||\partial_y(Q_{11},\varepsilon Q_{12})_\psi||_{\tilde L^2_{t,\eta'(t)}(B^\frac12)}||\varepsilon\partial_x u_\psi||_{\tilde L^2_{t,\eta'(t)}(B^\frac12)}.
\end{equation}
\end{prop}
\begin{proof}
From Poincar\'e inequality and $\partial_yv=-\partial_xu$, we have
\begin{equation*}
\begin{aligned}
&\int_0^t{(\Delta^h_k(T^h_{\partial_yQ_{11}}v)_\psi, \Delta^h_k(Q_{11})_\psi)_{L^2}dt'}\\
&\lesssim \sum_{|k'-k|\le 4}\int_0^t{||S_{k'-1}^h\partial_y(Q_{11})_\psi(t')||_{L^\infty}||\Delta^h_{k'}v_\psi(t')||_{L^2} ||\Delta^h_{k}(Q_{11})_\psi(t')||_{L^2}dt'}\\
&\lesssim\sum_{|k'-k|\le 4}\int_0^t{||\partial_y( Q_{11})_\psi(t')||_{B^\frac 12}||\Delta^h_{k'}v_\psi(t')||_{L^2} ||\Delta^h_{k}(Q_{11})_\psi(t')||_{L^2}dt'}\\
&\lesssim \sum_{|k'-k|\le 4}\left(\int_0^t{||\Delta^h_{k'}(\partial_xu)_\psi(t')||^2_{L^2} \eta'(t)dt'}\right)^\frac12\left(\int_0^t{||\Delta^h_{k}\partial_y(Q_{11})_\psi(t')||^2_{L^2}\eta'(t)dt'}\right)^\frac12\\
&\lesssim d_k^22^{-k}||\partial_xu_\psi||_{\tilde L^2_{t,\eta'(t)}(B^\frac12)}||\partial_y (Q_{11})_\psi||_{\tilde L^2_{t,\eta'(t)}(B^\frac12)}
\end{aligned}
\end{equation*}
Next,
\begin{equation*}
\begin{aligned}
&\int_0^t{(\Delta^h_k(T^h_v\partial_y Q_{11})_\psi, \Delta^h_k(Q_{11})_\psi)_{L^2}dt'}\\&\lesssim\sum_{|k'-k|\le 4}\int_0^t{||S_{k'-1}^hv_\psi(t')||_{L^\infty}|||\Delta^h_k(\partial_yQ_{11})_\psi(t')||_{L
^2} ||\Delta^h_{k}(Q_{11})_\psi(t')||_{L^2}dt'}\\
&\lesssim\sum_{|k'-k|\le 4}2^{-\frac{k'}{2}}\int_0^t{||S_{k'-1}^hv_\psi(t')||_{L^\infty}||(\partial_yQ_{11})_\psi(t')||_{B^\frac 12} ||\Delta^h_{k}(Q_{11})_\psi(t')||_{L^2}dt'}\\
&\lesssim \sum_{|k'-k|\le 4}2^{-\frac{k'}{2}}\left(\int_0^t{||S_{k'-1}^hv_\psi(t')||^2_{L^\infty} \eta'(t)dt'}\right)^\frac12\left(\int_0^t{||\Delta^h_{k}(\partial_y Q_{11})_\psi(t')||^2_{L^2}\eta'(t)dt'}\right)^\frac12\\
\end{aligned}
\end{equation*}
Notice that for any $s\in(0,\frac 12]$,
\begin{equation*}
\begin{aligned}
&\left(\int_0^t{||S_{k'-1}^hv_\psi(t')||^2_{L^\infty} \eta'(t)dt'}\right)^\frac12\\&\lesssim \sum_{l\le k'-2}2^\frac l2 \left(\int_0^t{||\Delta_{l}^hv_\psi(t')||^2_{L^2} \eta'(t)dt'}\right)^\frac12\lesssim \sum_{l\le k'-2}2^\frac l2 \left(\int_0^t{||\Delta_{l}^h\partial_xu_\psi(t')||^2_{L^2} \eta'(t)dt'}\right)^\frac12\\&\lesssim  \sum_{l\le k'-2}d_l2^{l(\frac 12-s)}||\partial_xu_\psi||_{\tilde L^2_{t,\eta'(t)}(B^s)} \lesssim  d_k2^{k(\frac 12-s)}||\partial_xu_\psi||_{\tilde L^2_{t,\eta'(t)}(B^s)}
\end{aligned}
\end{equation*}
hence
\[
\int_0^t{(\Delta^h_k(T^h_v\partial_y Q_{11})_\psi, \Delta^h_k(Q_{11})_\psi)_{L^2}dt'}\lesssim  d_k^22^{-k}||\partial_xu_\psi||_{\tilde L^2_{t,\eta'(t)}(B^\frac12)}||(\partial_yQ_{11})_\psi||_{\tilde L^2_{t,\eta'(t)}(B^\frac12)}
\]
For the third term, we have
\begin{equation*}
\begin{aligned}
&\int_0^t{(\Delta^h_k(R^h( v,\partial_y Q_{11}))_\psi, \Delta^h_k(Q_{11})_\psi)_{L^2}dt'}\\
&\lesssim 2^\frac k2\sum_{k'\ge k-3}\int_0^t{||\Delta^h_{k'}v_\psi(t')||_{L^2}||\tilde\Delta^h_{k'}(\partial_y Q_{11})_\psi(t')||_{L^2} ||\Delta^h_{k}(Q_{11})_\psi(t')||_{L^2}dt'}\\
&\lesssim 2^\frac{k}{2}\sum_{k'\ge k-3}2^{-\frac{k'}{2}}\int_0^t{||(\partial_yQ_{11})_\psi(t')||_{B^\frac 12}||\Delta^h_{k'}v_\psi(t')||_{L^2} ||\Delta^h_{k}(Q_{11})_\psi(t')||_{L^2}dt'}\\
&\lesssim 2^\frac{k}{2}\sum_{k'\ge k-3}2^{-\frac{k'}{2}}\int_0^t{||(\partial_yQ_{11})_\psi(t')||_{B^\frac 12}||\Delta^h_{k'}\partial_xu_\psi(t')||_{L^2} ||\Delta^h_{k}(\partial_yQ_{11})_\psi(t')||_{L^2}dt'}\\
&\lesssim2^\frac{k}{2}\sum_{k'\ge k-3}2^{-\frac{k'}{2}}\left(\int_0^t{||\Delta^h_{k'}( \partial_yQ_{11})_\psi(t')||^2_{L^2} \eta'(t)dt'}\right)^\frac12\left(\int_0^t{||\Delta^h_{k}\partial_xu_\psi(t')||^2_{L^2}\eta'(t)dt'}\right)^\frac12\\
&\lesssim d_k2^{-k}||\partial_y(Q_{11})_\psi||_{\tilde L^2_{t,\eta'(t)}(B^\frac12)}||\partial_x u_\psi||_{\tilde L^2_{t,\eta'(t)}(B^\frac12)}\left(\sum_{k'\ge k-3}d_{k'}2^{(k-k')}\right)\\
&\lesssim d_k^22^{-k}||\partial_y(Q_{11})_\psi||_{\tilde L^2_{t,\eta'(t)}(B^\frac12)}||\partial_x u_\psi||_{\tilde L^2_{t,\eta'(t)}(B^\frac12)}
\end{aligned}
\end{equation*}
and the proposition is proved.
\end{proof}
\begin{prop}
\begin{equation}
\label{estib7}
B_3\lesssim d_k^22^{-k}\left(||e^{\mathcal Rt'}(\varepsilon\partial_xu)_\psi||^2_{\tilde L^2_{t,\eta'(t)}(B^\frac12)}+||e^{\mathcal Rt'}(\partial_y(Q_{11}, \varepsilon Q_{12})_\psi||_{\tilde L^2_{t,\eta'(t)}(B^\frac12)}^2\right)
\end{equation}
\end{prop}
\begin{proof}
The terms related to $\partial_yu$ can be estimated by using Lemma \ref{lem34}. Next, recall from (3.18) of \cite{PZZ} that because $v(t,x,y)=-\int_0^y{\partial_xu(t,x,y')dy'}$, we have
\begin{equation}
\label{extrav1}
||\Delta^h_kv_\psi(t)||_{L^\infty}\lesssim 2^{\frac {3k}2}||\Delta^h_ku_\psi(t)||_{L^2}
\end{equation}
Notice that
\begin{equation*}
\begin{aligned}
&\varepsilon^3\int_0^t{(\Delta^h_k(T^h_{\partial_xv}Q_{12})_\psi, \Delta^h_k(Q_{11})_\psi)_{L^2}dt'}\\
&\lesssim\varepsilon ^3\sum_{|k'-k|\le 4}\int_0^t{||S_{k'-1}^h\partial_xv_\psi(t')||_{L^\infty}||\Delta^h_{k'}(Q_{12})_\psi(t')||_{L^2} ||\Delta^h_{k}(Q_{11})_\psi(t')||_{L^2}dt'}\\
&\lesssim\varepsilon\sum_{|k'-k|\le 4}2^{-\frac{3k'}2}\int_0^t{||S_{k'-1}^h\partial_xv_\psi(t')||_{L^\infty}||\Delta^h_{k'}(\varepsilon^2\partial_xQ_{12})_\psi(t')||_{B^\frac12} ||\Delta^h_{k}(\partial_yQ_{11})_\psi(t')||_{L^2}dt'}\\
&\lesssim \sum_{|k'-k|\le 4}2^{-\frac{3k'}2}\varepsilon\left(\int_0^t{||S_{k'-1}^h\partial_xv_\psi(t')||_{L^\infty}^2 \eta'(t)dt'}\right)^\frac12\left(\int_0^t{||\Delta^h_{k}\partial_y(Q_{11})_\psi(t')||^2_{L^2}\eta'(t)dt'}\right)^\frac12\\
\end{aligned}
\end{equation*}
and we deduce from \eqref{extrav1} that
\begin{equation*}
\begin{aligned}
\left(\int_0^t{||S_{k'-1}^h\varepsilon\partial_xv_\psi(t')||_{L^\infty}^2 \eta'(t)dt'}\right)^\frac12&\lesssim\sum_{l\le k'-2}2^\frac{3l}{2}\left(\int_0^t{||\Delta^h_l\varepsilon\partial_xu_\psi(t')||_{L^2}^2 \eta'(t)dt'}\right)^\frac12\\
&\lesssim\sum_{l\le k'-2}d_l2^{l}||\varepsilon\partial_xu_\psi||_{\tilde L^2_{t,\eta'(t)}(B^\frac 12)}\lesssim d_{k'}2^{k'}||\varepsilon\partial_xu_\psi||_{\tilde L^2_{t,\eta'(t)}(B^\frac 12)}
\end{aligned}
\end{equation*}
so
\[
\varepsilon^3\int_0^t{(\Delta^h_k(T^h_{\partial_xv}Q_{12})_\psi, \Delta^h_k(Q_{11})_\psi)_{L^2}dt'}\lesssim d_k^22^{-k}||\varepsilon\partial_xu_\psi||_{\tilde L^2_{t,\eta'(t)}(B^\frac 12)}||\partial_y (Q_{11})_\psi||_{\tilde L^2_{t,\eta'(t)}(B^\frac12)}
\]
Similarly, we have
\begin{equation*}
\begin{aligned}
&\varepsilon^3\int_0^t{(\Delta^h_k(T^h_{Q_{12}}\partial_xv)_\psi, \Delta^h_k(Q_{11})_\psi)_{L^2}dt'}\\&\lesssim\varepsilon^3\sum_{|k'-k|\le 4}\int_0^t{||S_{k'-1}^h( Q_{12})_\psi(t')||_{L^\infty}|||\Delta^h_{k'}(\partial_xv)_\psi(t')||_{L
^2} ||\Delta^h_{k}(Q_{11})_\psi(t')||_{L^2}dt'}\\
&\lesssim\sum_{|k'-k|\le 4}\int_0^t{||(\varepsilon^2 Q_{12})_\psi(t')||_{B^\frac 32}||(\varepsilon\Delta^h_{k'}\partial_xu)_\psi(t')||_{L^2} ||\Delta^h_{k}(\partial_y Q_{11})_\psi(t')||_{L^2}dt'}\\
&\lesssim \sum_{|k'-k|\le 4}\left(\int_0^t{||(\varepsilon\partial_xu)_\psi(t')||_{L^2}^2 \eta'(t)dt'}\right)^\frac12\left(\int_0^t{||\Delta^h_{k}(\partial_y Q_{11})_\psi(t')||^2_{L^2}\eta'(t)dt'}\right)^\frac12\\
&\lesssim d_k^22^{-k}||\varepsilon\partial_xu_\psi||_{\tilde L^2_{t,\eta'(t)}(B^\frac12)}||\partial_y (Q_{11})_\psi||_{\tilde L^2_{t,\eta'(t)}(B^\frac12)}
\end{aligned}
\end{equation*}
and
\begin{equation*}
\begin{aligned}
&\varepsilon^3\int_0^t{(\Delta^h_k(R^h( \partial_xv, Q_{12}))_\psi, \Delta^h_k(Q_{11})_\psi)_{L^2}dt'}\\
&\lesssim \varepsilon^32^\frac k2\sum_{k'\ge k-3}\int_0^t{||\Delta^h_{k'}\partial_xv_\psi(t')||_{L^2}||\tilde\Delta^h_{k'}( Q_{12})_\psi(t')||_{L^2} ||\Delta^h_{k}(\partial_yQ_{11})_\psi(t')||_{L^2}dt'}\\
&\lesssim 2^\frac{k}{2}\sum_{k'\ge k-3}2^{-\frac{k'}{2}}\int_0^t{||(\varepsilon\partial_xu)_\psi(t')||_{L^2}||\Delta^h_{k'}(\varepsilon^2\partial_xQ_{12})_\psi(t')||_{B^\frac 12} ||\Delta^h_{k}(\partial_yQ_{11})_\psi(t')||_{L^2}dt'}\\
&\lesssim2^\frac{k}{2} \sum_{k'\ge k-3}2^{-\frac{k'}{2}}\left(\int_0^t{||\Delta^h_{k'}\varepsilon\partial_xu_\psi(t')||^2_{L^2} \eta'(t)dt'}\right)^\frac12\left(\int_0^t{||\Delta^h_{k}\partial_y(Q_{11})_\psi(t')||^2_{L^2}\eta'(t)dt'}\right)^\frac12\\
&\lesssim d_k2^{-k}||\varepsilon\partial_xu_\psi||_{\tilde L^2_{t,\eta'(t)}(B^\frac12)}||\partial_y(Q_{11})_\psi||_{\tilde L^2_{t,\eta'(t)}(B^\frac12)}\left(\sum_{k'\ge k-3}d_{k'}2^{(k-k')}\right)\\
&\lesssim d_k^22^{-k}||\varepsilon\partial_xu_\psi||_{\tilde L^2_{t,\eta'(t)}(B^\frac12)}||\partial_y(Q_{11})_\psi||_{\tilde L^2_{t,\eta'(t)}(B^\frac12)}
\end{aligned}
\end{equation*}
which together with Lemma \ref{lem34} provides the result.
\end{proof}

\begin{rmk}
The term $B_4$ can be estimated similarly as $B_3$ by exchanging the order of $Q_{11}$ and $Q_{12}$.
\end{rmk}

\subsubsection{Estimates of $A_5, A_6, A_7, A_8$} We prove the following propositions:
\begin{prop}
\label{propa5} We have
\begin{equation}
\label{estia7}
A_5+ A_8\lesssim  \varepsilon d_k^22^{-k}||e^{\mathcal R't'} (\varepsilon\partial_x,\partial_y)(Q_{11}, \varepsilon Q_{12})_\psi||_{\tilde L^2_{t,\eta'(t)}(B^\frac12)}||e^{\mathcal R't'}\varepsilon\partial_x u_\psi||_{\tilde L^2_{t,\eta'(t)}(B^\frac12)}
\end{equation}
\end{prop}
\begin{proof}
Notice that
\begin{equation*}
\begin{aligned}
&\varepsilon^3\int_0^t{(\Delta^h_k(T^h_{\partial_xQ_{11}}\partial_xQ_{11})_\psi, \Delta^h_k\partial_xu_\psi)_{L^2}dt'}\\
&\lesssim\varepsilon^3\sum_{|k'-k|\le 4}\int_0^t{||S_{k'-1}^h\partial_x(Q_{11})_\psi(t')||_{L^\infty}||\Delta^h_{k'}(\partial_xQ_{11})_\psi(t')||_{L^2} ||\Delta^h_{k}\partial_xu_\psi(t')||_{L^2}dt'}\\
&\lesssim\sum_{|k'-k|\le 4}\int_0^t{||(\varepsilon \partial_xQ_{11})_\psi(t')||_{B^\frac 12}||\Delta^h_{k'}(
\varepsilon\partial_xQ_{11})_\psi(t')||_{L^2} ||\Delta^h_{k}\varepsilon\partial_xu_\psi(t')||_{L^2}dt'}\\
&\lesssim \sum_{|k'-k|\le 4}\left(\int_0^t{||\Delta^h_{k'}(\varepsilon\partial_xQ_{11})_\psi(t')||^2_{L^2} \eta'(t)dt'}\right)^\frac12\left(\int_0^t{||\Delta^h_{k}\varepsilon\partial_xu_\psi(t')||^2_{L^2}\eta'(t)dt'}\right)^\frac12\\
&\lesssim d_k^22^{-k}||\varepsilon\partial_x(Q_{11})_\psi||_{\tilde L^2_{t,\eta'(t)}(B^\frac12)}||\varepsilon\partial_x u_\psi||_{\tilde L^2_{t,\eta'(t)}(B^\frac12)}
\end{aligned}
\end{equation*}
and
\begin{equation*}
\begin{aligned}
&\varepsilon^3\int_0^t{(\Delta^h_k(R^h({Q_{11}},\partial_y Q_{12}))_\psi, \Delta^h_k\partial_x^2u_\psi)_{L^2}dt'}\\
&\lesssim\varepsilon^32^\frac k2\sum_{k'\ge k-3}\int_0^t{||\Delta^h_{k'}(\partial_xQ_{11})_\psi(t')||_{L^2}||\tilde\Delta^h_{k'}(\partial_x Q_{11})_\psi(t')||_{L^2} ||\Delta^h_{k}\partial_xu_\psi(t')||_{L^2}dt'}\\
&\lesssim2^\frac{k}{2}\sum_{|k'-k|\le 4}2^{-\frac{k'}{2}}\int_0^t{||(\varepsilon\partial_xQ_{11})_\psi(t')||_{B^\frac 12}||\Delta^h_{k'}(\varepsilon\partial_xQ_{11})_\psi(t')||_{L^2} ||\Delta^h_{k}\varepsilon\partial_xu_\psi(t')||_{L^2}dt'}\\
&\lesssim 2^\frac{k}{2}\sum_{|k'-k|\le 4}2^\frac{-k'}{2}\left(\int_0^t{||\Delta^h_{k'}(\varepsilon \partial_xQ_{11})_\psi(t')||^2_{L^2} \eta'(t)dt'}\right)^\frac12\left(\int_0^t{||\Delta^h_{k}\varepsilon\partial_xu_\psi(t')||^2_{L^2}\eta'(t)dt'}\right)^\frac12\\
&\lesssim d_k2^{-k}||\partial_x(Q_{11})_\psi||_{\tilde L^2_{t,\eta'(t)}(B^\frac12)}||\varepsilon\partial_x u_\psi||_{\tilde L^2_{t,\eta'(t)}(B^\frac12)}\left(\sum_{k'\ge k-3}d_{k'}2^{2(k-k')}\right)\\
&\lesssim d_k^22^{-k}||\partial_x(Q_{11})_\psi||_{\tilde L^2_{t,\eta'(t)}(B^\frac12)}||\varepsilon\partial_x u_\psi||_{\tilde L^2_{t,\eta'(t)}(B^\frac12)}
\end{aligned}
\end{equation*}
so we have proved that \[A_5\lesssim  \varepsilon d_k^22^{-k}||e^{\mathcal R't'} \varepsilon\partial_x(Q_{11}, \varepsilon Q_{12})_\psi||_{\tilde L^2_{t,\eta'(t)}(B^\frac12)}||e^{\mathcal R't'}\varepsilon\partial_x u_\psi||_{\tilde L^2_{t,\eta'(t)}(B^\frac12)}\]
The estimate of $A_8$ is similar, just exchanging $\varepsilon \partial_x$ to $\partial_y$ and using  that $\partial_yv=-\partial_xu$.\end{proof}

\begin{prop}
\begin{equation}
\label{equationa9}
\begin{aligned}
&A_7\lesssim \varepsilon d_k^22^{-k}\left(||\varepsilon\partial_x(Q_{11}, \varepsilon Q_{12})_\psi||_{\tilde L^2_{t,\eta'(t)}(B^\frac12)}+||\partial_y(Q_{11}, \varepsilon Q_{12})_\psi||_{\tilde L^2_{t,\eta'(t)}(B^\frac12)}\right)||(\varepsilon\partial_xu, \varepsilon^2\partial_x v)_\psi||_{\tilde L^2_{t,\eta'(t)}(B^\frac12)}
\end{aligned}
\end{equation}
\end{prop}
\begin{proof}We write
\[
\Delta_{\varepsilon}Q_{12}Q_{11}-\Delta_{\varepsilon}Q_{11}Q_{12}=\varepsilon^2\partial_x^2Q_{12}Q_{11}-\varepsilon^2\partial_x^2Q_{11}Q_{12}+\partial_y^2Q_{12}Q_{11}-\partial_y^2Q_{11}Q_{12}
\]
from integration by parts,
\[
(\Delta^h_k (\partial_y^2Q_{12}Q_{11})_\psi, \Delta^h_k \partial_x v_\psi)_{L^2}+(\Delta^h_k (\partial_yQ_{12}\partial_yQ_{11})_\psi, \Delta^h_k \partial_x v_\psi)_{L^2}=-(\Delta^h_k (\partial_yQ_{12}Q_{11})_\psi, \Delta^h_k \partial_x\partial_y v_\psi)_{L^2}
\]
so we have (recall that $\partial_yv=-\partial_xu$)
\[
(\Delta^h_k (\partial_y^2Q_{12}Q_{11}-\partial_y^2Q_{11}Q_{12})_\psi, \Delta^h_k \partial_x v_\psi)_{L^2}=(\Delta^h_k (\partial_yQ_{12}Q_{11})_\psi, \Delta^h_k \partial_x^2 u_\psi)_{L^2}-(\Delta^h_k (\partial_yQ_{11}Q_{12})_\psi, \Delta^h_k \partial_x^2 u_\psi)_{L^2}
\]
Notice that
\begin{equation*}
\begin{aligned}
&\varepsilon^3\int_0^t{(\Delta^h_k(T^h_{\partial_yQ_{12}}Q_{11})_\psi, \Delta^h_k\partial_x^2u_\psi)_{L^2}dt'}\\
&\lesssim\varepsilon^3\sum_{|k'-k|\le 4}\int_0^t{||S_{k'-1}^h\partial_y(Q_{12})_\psi(t')||_{L^\infty}||\Delta^h_{k'}(Q_{11})_\psi(t')||_{L^2} ||\Delta^h_{k}\partial_x^2u_\psi(t')||_{L^2}dt'}\\
&\lesssim\varepsilon^2\sum_{|k'-k|\le 4}\int_0^t{||\partial_y(\varepsilon Q_{12})_\psi(t')||_{B^\frac 12}||\Delta^h_{k'}(Q_{11})_\psi(t')||_{L^2} ||\Delta^h_{k}\partial_x^2u_\psi(t')||_{L^2}dt'}\\
&\lesssim \sum_{|k'-k|\le 4}\left(\int_0^t{||\Delta^h_{k'}(\varepsilon\partial_xQ_{11})_\psi(t')||^2_{L^2} \eta'(t)dt'}\right)^\frac12\left(\int_0^t{||\Delta^h_{k}\varepsilon\partial_xu_\psi(t')||^2_{L^2}\eta'(t)dt'}\right)^\frac12\\
&\lesssim d_k^22^{-k}||\varepsilon\partial_x(Q_{11})_\psi||_{\tilde L^2_{t,\eta'(t)}(B^\frac12)}||\varepsilon\partial_x u_\psi||_{\tilde L^2_{t,\eta'(t)}(B^\frac12)}
\end{aligned}
\end{equation*}
and similarly,
\begin{equation*}
\begin{aligned}
&\varepsilon^3\int_0^t{(\Delta^h_k(T^h_{Q_{11}}\partial_y Q_{12})_\psi, \Delta^h_k\partial_x^2u_\psi)_{L^2}dt'}\\
&\lesssim\varepsilon^3\sum_{|k'-k|\le 4}\int_0^t{||S_{k'-1}^h(Q_{11})_\psi(t')||_{L^\infty}||\Delta^h_{k'}(\partial_y Q_{12})_\psi(t')||_{L^2} ||\Delta^h_{k}\partial_x^2u_\psi(t')||_{L^2}dt'}\\
&\lesssim\varepsilon^2\sum_{|k'-k|\le 4}2^{-k'}\int_0^t{||(\varepsilon Q_{11})_\psi(t')||_{B^\frac 32}||\Delta^h_{k'}(\partial_yQ_{12})_\psi(t')||_{L^2} ||\Delta^h_{k}\partial_x^2u_\psi(t')||_{L^2}dt'}\\
&\lesssim \sum_{|k'-k|\le 4}\left(\int_0^t{||\Delta^h_{k'}(\partial_yQ_{12})_\psi(t')||^2_{L^2} \eta'(t)dt'}\right)^\frac12\left(\int_0^t{||\Delta^h_{k}\varepsilon\partial_xu_\psi(t')||^2_{L^2}\eta'(t)dt'}\right)^\frac12\\
&\lesssim d_k^22^{-k}||\partial_y(Q_{12})_\psi||_{\tilde L^2_{t,\eta'(t)}(B^\frac12)}||\varepsilon\partial_x u_\psi||_{\tilde L^2_{t,\eta'(t)}(B^\frac12)}
\end{aligned}
\end{equation*}
Moreover,
\begin{equation*}
\begin{aligned}
&\varepsilon^2\int_0^t{(\Delta^h_k(R^h({Q_{11}},\partial_y Q_{12}))_\psi, \Delta^h_k\partial_x^2u_\psi)_{L^2}dt'}\\
&\lesssim\varepsilon^22^\frac k2\sum_{k'\ge k-3}\int_0^t{||\Delta^h_{k'}(Q_{11})_\psi(t')||_{L^2}||\tilde\Delta^h_{k'}(\partial_y Q_{12})_\psi(t')||_{L^2} ||\Delta^h_{k}\partial_x^2u_\psi(t')||_{L^2}dt'}\\
&\lesssim\varepsilon^22^\frac{3k}{2}\sum_{|k'-k|\le 4}2^{-\frac{3k'}{2}}\int_0^t{||(\partial_yQ_{12})_\psi(t')||_{B^\frac 12}||\Delta^h_{k'}(\partial_xQ_{11})_\psi(t')||_{L^2} ||\Delta^h_{k}\partial_xu_\psi(t')||_{L^2}dt'}\\
&\lesssim2^\frac{3k}{2} \sum_{|k'-k|\le 4}2^\frac{-3k'}{2}\left(\int_0^t{||\Delta^h_{k'}(\varepsilon \partial_xQ_{11})_\psi(t')||^2_{L^2} \eta'(t)dt'}\right)^\frac12\left(\int_0^t{||\Delta^h_{k}\varepsilon\partial_xu_\psi(t')||^2_{L^2}\eta'(t)dt'}\right)^\frac12\\
&\lesssim d_k2^{-k}||\partial_x(Q_{11})_\psi||_{\tilde L^2_{t,\eta'(t)}(B^\frac12)}||\varepsilon\partial_x u_\psi||_{\tilde L^2_{t,\eta'(t)}(B^\frac12)}\left(\sum_{k'\ge k-3}d_{k'}2^{2(k-k')}\right)\\
&\lesssim d_k^22^{-k}||\partial_x(Q_{11})_\psi||_{\tilde L^2_{t,\eta'(t)}(B^\frac12)}||\varepsilon\partial_x u_\psi||_{\tilde L^2_{t,\eta'(t)}(B^\frac12)}
\end{aligned}
\end{equation*}
so after combing these three terms, we have
\begin{equation*}
\begin{aligned}
&\varepsilon^3\int_0^t{(\Delta^h_k(\partial_yQ_{12}Q_{11})_\psi, \Delta^h_k\partial_x^2u_\psi)_{L^2}dt'}\\&\lesssim d_k^22^{-k}\left(||\varepsilon\partial_x(Q_{11})_\psi||_{\tilde L^2_{t,\eta'(t)}(B^\frac12)}+||\partial_y(\varepsilon Q_{12})_\psi||_{\tilde L^2_{t,\eta'(t)}(B^\frac12)}\right)||\varepsilon\partial_x u_\psi||_{\tilde L^2_{t,\eta'(t)}(B^\frac12)}
\end{aligned}
\end{equation*}
Next,
\[
(\Delta^h_k (\partial_x^2Q_{12}Q_{11}-\partial_x^2Q_{11}Q_{12})_\psi, \Delta^h_k \partial_x v_\psi)_{L^2}=(\Delta^h_k (\partial_xQ_{12}Q_{11})_\psi, \Delta^h_k \partial_x^2 v_\psi)_{L^2}-(\Delta^h_k (\partial_xQ_{11}Q_{12})_\psi, \Delta^h_k \partial_x^2 v_\psi)_{L^2}
\]
Similarly as the proof above (just change $\partial_yQ_{12}$ to $\varepsilon\partial_xQ_{12}$ and $\partial_xu$ to $\varepsilon\partial_xv$), we have
\begin{equation*}
\begin{aligned}
&\varepsilon^5\int_0^t{(\Delta^h_k(\partial_xQ_{12}Q_{11})_\psi, \Delta^h_k\partial_x^2v_\psi)_{L^2}dt'}\\&\lesssim d_k^22^{-k}\left(||\varepsilon\partial_x(Q_{11})_\psi||_{\tilde L^2_{t,\eta'(t)}(B^\frac12)}+||\partial_x(\varepsilon^2 Q_{12})_\psi||_{\tilde L^2_{t,\eta'(t)}(B^\frac12)}\right)||\varepsilon^2\partial_x v_\psi||_{\tilde L^2_{t,\eta'(t)}(B^\frac12)}
\end{aligned}
\end{equation*}
and then
\begin{equation*}
\begin{aligned}
&\varepsilon(\Delta^h_k(R^\varepsilon_{12,2})_\psi, \Delta^h_k\partial_xv_\psi)_{L^2}\lesssim d_k^22^{-k}||\varepsilon\partial_x(Q_{11},\varepsilon Q_{12})_\psi||_{\tilde L^2_{t,\eta'(t)}(B^\frac12)}||\varepsilon^2\partial_x v_\psi||_{\tilde L^2_{t,\eta'(t)}(B^\frac12)}
\\
&+d_k^22^{-k}\left(||\varepsilon\partial_x(Q_{11}, \varepsilon Q_{12})_\psi||_{\tilde L^2_{t,\eta'(t)}(B^\frac12)}
+||\partial_y(Q_{11}, \varepsilon Q_{12})_\psi||_{\tilde L^2_{t,\eta'(t)}(B^\frac12)}\right)||\varepsilon\partial_x u_\psi||_{\tilde L^2_{t,\eta'(t)}(B^\frac12)}\end{aligned}
\end{equation*}
Similarly, we have
\begin{equation*}
\begin{aligned}
\varepsilon(\Delta^h_k(R^\varepsilon_{12,1})_\psi,\,\,& \Delta^h_k\partial_xv_\psi)_{L^2}\\&\lesssim d_k^22^{-k}\left(||\varepsilon\partial_x(Q_{11}, \varepsilon Q_{12})_\psi||_{\tilde L^2_{t,\eta'(t)}(B^\frac12)}+||\partial_y(Q_{11}, \varepsilon Q_{12})_\psi||_{\tilde L^2_{t,\eta'(t)}(B^\frac12)}\right)||\varepsilon^2\partial_x v_\psi||_{\tilde L^2_{t,\eta'(t)}(B^\frac12)}\\
\end{aligned}
\end{equation*}
and the proposition is proved by gathering the estimates above.
\end{proof}

\begin{prop}
\label{propa8}
\begin{equation}
\begin{aligned}
\label{estia8}
&A_6\lesssim \varepsilon^3d_k^22^{-k}||\partial_x(Q_{11},\varepsilon Q_{12})_\psi||_{\tilde L^2_{t,\eta'(t)}(B^{\frac 12})}||\partial_y(Q_{11},\varepsilon Q_{12})_\psi||_{\tilde L^2_{t,\eta'(t)}(B^{\frac 12})}\\
&\lesssim \varepsilon^3d_k^22^{-k}||(Q_{11})_\psi||_{\tilde L^2_{t,\eta'(t)}(B^{\frac 12})}||(\Delta_\varepsilon Q_{12})_\psi||_{\tilde L^2_{t,\eta'(t)}(B^{\frac 12})}+\varepsilon^3d_k^22^{-k}||(Q_{12})_\psi||_{\tilde L^2_{t,\eta'(t)}(B^{\frac 12})}||(\Delta_\varepsilon Q_{11})_\psi||_{\tilde L^2_{t,\eta'(t)}(B^{\frac 12})}.\\
\end{aligned}
\end{equation}
\end{prop}
\begin{proof}
The proof is similar as Lemma \ref{lem32}. We skip it for simplicity.
\end{proof}

From the estimate of $A_6$, we notice that there exists an extra term $\Delta_\varepsilon (Q_{11}, \varepsilon Q_{12})$ (more precisely, it is $\partial^2_y (Q_{11}, \varepsilon Q_{12})$ ), which could be estimated from the calculations above. To solve this problem, we come back to the equation \eqref{q11v2}. It has the term $\Delta_\varepsilon Q$ on the right, so we apply $\Delta^h_k $ to it and take the $L^2$ inner product with $-\varepsilon^2\Delta^h_k(\Delta_\varepsilon (Q_{11}, \varepsilon Q_{22}))_\psi$ to provide the bound related to $\Delta_\varepsilon Q$. The reason that we multiply by an extra $\varepsilon^2$ is to prove the convergence to the hydrostatic system, to be done later. We have
\begin{equation}
\label{deltaq3}
\begin{aligned}
&\frac12\frac d{dt}||(\varepsilon^2 \partial_x, \varepsilon\partial_y)\Delta^h_k(Q_{11},\varepsilon Q_{22})_\psi||^2_{L^2}+a||(\varepsilon^2 \partial_x, \varepsilon\partial_y)\Delta^h_k(Q_{11},\varepsilon Q_{22})_\psi||^2_{L^2}\\&+\lambda \eta'(t')(|D_x|\Delta^h_k(\varepsilon^2 \partial_x, \varepsilon\partial_y)(Q_{11},\varepsilon Q_{22})_\psi,\Delta^h_k(\varepsilon^2 \partial_x, \varepsilon\partial_y)(Q_{11},\varepsilon Q_{22})_\psi)_{L^2} +||\varepsilon \Delta^h_k(\Delta_\varepsilon Q_{11},\varepsilon \Delta_\varepsilon Q_{12})_\psi||_{L^2}^2\\
&=\varepsilon^2(\Delta^h_k((U+\varepsilon u)\partial_x Q_{11})_\psi, \Delta^h_k(\Delta_\varepsilon Q_{11})_\psi)_{L^2}+\varepsilon^4(\Delta^h_k((U+\varepsilon u)\partial_x Q_{12})_\psi, \Delta^h_k(\Delta_\varepsilon Q_{12})_\psi)_{L^2}\\
&+\varepsilon^3(\Delta^h_k((v\partial_y Q_{11})_\psi, \Delta^h_k(\Delta_\varepsilon Q_{11})_\psi)_{L^2}+\varepsilon^5(\Delta^h_k((v\partial_y Q_{12})_\psi, \Delta^h_k(\Delta_\varepsilon Q_{12})_\psi)_{L^2}\\
&+ \varepsilon^2(\Delta^h_k((\partial_y(U+\varepsilon u) Q_{12})_\psi, \Delta^h_k(\Delta_\varepsilon Q_{11})_\psi)_{L^2}+ \varepsilon^2(\Delta^h_k((\partial_y(U+\varepsilon u) Q_{11})_\psi, e^{\mathcal R't'}\Delta^h_k(\Delta_\varepsilon Q_{12})_\psi)_{L^2}\\
&- \varepsilon^5(\Delta^h_k((\partial_xv Q_{12})_\psi, \Delta^h_k(\Delta_\varepsilon Q_{11})_\psi)_{L^2}- \varepsilon^5(\Delta^h_k((\partial_xvQ_{11})_\psi, \Delta^h_k(\Delta_\varepsilon Q_{12})_\psi)_{L^2}\\
&+2c\varepsilon^2(\Delta^h_k(Q_{11}^3+2\varepsilon^2Q_{11}Q_{12}^2)_{\psi}, \Delta^h_k(\Delta_\varepsilon Q_{11})_\psi)_{L^2}+2c\varepsilon^4(\Delta^h_k(2\varepsilon^2Q_{12}^3+Q_{12}Q_{11}^2)_{\psi}, \Delta^h_k(\Delta_\varepsilon Q_{12})_\psi)_{L^2}\\
\end{aligned}
\end{equation}
Multiply \eqref{deltaq3} by $e^{2\mathcal Rt}$ and integrate over $[0,t]$, to get
\begin{equation}
\label{deltaq3}
\begin{aligned}
&\frac12||e^{\mathcal Rt'}(\varepsilon^2 \partial_x, \varepsilon\partial_y)\Delta^h_k(Q_{11},\varepsilon Q_{22})_\psi||^2_{L^\infty_t(L^2)}+a'\int_0^t{||e^{\mathcal Rt'}(\varepsilon^2 \partial_x, \varepsilon\partial_y)\Delta^h_k(Q_{11},\varepsilon Q_{22})_\psi(t')||^2_{L^2}dt'}\\
&+\lambda 2^k\int_0^t{\eta'(t)||e^{\mathcal Rt'}(\varepsilon^2 \partial_x, \varepsilon\partial_y)\Delta^h_k(Q_{11},\varepsilon Q_{22})_\psi(t')||^2_{L^2}dt'}+\int_0^t{||e^{\mathcal Rt'}\varepsilon \Delta^h_k(\Delta_\varepsilon Q_{11},\varepsilon \Delta_\varepsilon Q_{12})_\psi(t')||_{L^2}^2dt'}\\
&\le ||e^{a|D_x|}(\varepsilon^2\partial_x,\varepsilon\partial_y)\Delta^h_k((Q_{11})_0,\varepsilon (Q_{12})_0)||_{L^2}^2+\underbrace{\varepsilon^2\int_0^t{|e^{\mathcal Rt'}(\Delta^h_k((U+\varepsilon u)\partial_x Q_{11})_\psi, e^{\mathcal Rt'}\Delta^h_k(\Delta_\varepsilon Q_{11})_\psi)_{L^2}|dt'}}_{C_1}\\
&+\underbrace{\varepsilon^4\int_0^t{|e^{\mathcal Rt'}(\Delta^h_k((U+\varepsilon u)\partial_x Q_{12})_\psi,e^{\mathcal Rt'} \Delta^h_k(\Delta_\varepsilon Q_{12})_\psi)_{L^2}|dt'}}_{C_2}+\underbrace{\varepsilon^3\int_0^t{|(e^{\mathcal Rt'}\Delta^h_k((v\partial_y Q_{11})_\psi,e^{\mathcal Rt'} \Delta^h_k(\Delta_\varepsilon Q_{11})_\psi)_{L^2}|dt'}}_{C_3}\\
&+\underbrace{\varepsilon^5\int_0^t{|e^{\mathcal Rt'}(\Delta^h_k((v\partial_y Q_{12})_\psi,e^{\mathcal Rt'} \Delta^h_k(\Delta_\varepsilon Q_{12})_\psi)_{L^2}|dt'}}_{C_4}
+\underbrace{ \varepsilon^2\int_0^t{|(e^{\mathcal Rt'}\Delta^h_k((\partial_y(U+\varepsilon u) Q_{12})_\psi,e^{\mathcal Rt'} \Delta^h_k(\Delta_\varepsilon Q_{11})_\psi)_{L^2}|dt'}}_{C_5}\\
&+ \underbrace{\varepsilon^2\int_0^t{|(e^{\mathcal Rt'}\Delta^h_k((\partial_y(U+\varepsilon u) Q_{11})_\psi,e^{\mathcal Rt'} \Delta^h_k(\Delta_\varepsilon Q_{12})_\psi)_{L^2}|dt'}}_{C_6}+ \underbrace{\varepsilon^5\int_0^t{|(e^{\mathcal Rt'}\Delta^h_k((\partial_xv Q_{12})_\psi, e^{\mathcal Rt'}\Delta^h_k(\Delta_\varepsilon Q_{11})_\psi)_{L^2}|dt'}}_{C_7}\\
&+ \underbrace{\varepsilon^5\int_0^t{|(e^{\mathcal Rt'}\Delta^h_k((\partial_xv Q_{11})_\psi, e^{\mathcal Rt'}\Delta^h_k(\Delta_\varepsilon Q_{12})_\psi)_{L^2}|dt'}}_{C_8}+\underbrace{2c\varepsilon^2\int_0^t{|(e^{\mathcal Rt'}\Delta^h_k(Q_{11}^3+2\varepsilon^2Q_{11}Q_{12}^2)_{\psi}, e^{\mathcal Rt'}\Delta^h_k(\Delta_\varepsilon Q_{11})_\psi)_{L^2}|dt'}}_{C_9}\\
&+\underbrace{2c'\varepsilon^4\int_0^t{|(e^{\mathcal Rt'}\Delta^h_k(2\varepsilon^2Q_{12}^3+Q_{12}Q_{11}^2)_{\psi}, e^{\mathcal Rt'}\Delta^h_k(\Delta_\varepsilon Q_{12})_\psi)_{L^2}|dt'}}_{C_{10}}\\
\end{aligned}
\end{equation}

Similarly as  the estimate of $B_1$ and recalling that $U\lesssim \eta'(t)$, we have

\[
C_1+C_2\lesssim \varepsilon^2d_k^22^{-k}||e^{\mathcal Rt'}\partial_x(Q_{11},\varepsilon Q_{12})_\psi||_{\tilde L^2_{t,\eta'(t)}(B^\frac 12)}||e^{\mathcal Rt'}\Delta_\varepsilon(Q_{11},\varepsilon Q_{12})_\psi||_{\tilde L^2_{t,\eta'(t)}(B^\frac 12)}
\]
Similarly as the estimate of $B_2$, we get

\[
C_3+C_4\lesssim \varepsilon^3d_k^22^{-k}||e^{\mathcal Rt'}\partial_xu_\psi||_{\tilde L^2_{t,\eta'(t)}(B^\frac 12)}||e^{\mathcal Rt'}\Delta_\varepsilon(Q_{11},\varepsilon Q_{12})_\psi||_{\tilde L^2_{t,\eta'(t)}(B^\frac 12)}
\]
Similarly as the estimate of $B_3$, and recalling that $\partial_yU\lesssim \eta'(t)$, we have
\begin{equation*}
\begin{aligned}
C_5+C_6&\lesssim \varepsilon^2d_k^22^{-k}(||e^{\mathcal Rt'}(Q_{12})_\psi||_{\tilde L^2_{t,\eta'(t)}(B^\frac 12)}||e^{\mathcal Rt'}\Delta_\varepsilon (Q_{11})_\psi||_{\tilde L^2_{t,\eta'(t)}(B^\frac 12)})\\
&+\varepsilon^2d_k^22^{-k}(||e^{\mathcal Rt'}(Q_{11})_\psi||_{\tilde L^2_{t,\eta'(t)}(B^\frac 12)}||e^{\mathcal Rt'}\Delta_\varepsilon (Q_{12})_\psi||_{\tilde L^2_{t,\eta'(t)}(B^\frac 12)})
\end{aligned}
\end{equation*}
and
\[
C_7+C_8\lesssim \varepsilon^3d_k^22^{-k}||e^{\mathcal Rt'}\partial_xu_\psi||_{\tilde L^2_{t,\eta'(t)}(B^\frac 12)}||e^{\mathcal Rt'}\Delta_\varepsilon(Q_{11},\varepsilon Q_{12})_\psi||_{\tilde L^2_{t,\eta'(t)}(B^\frac 12)}
\]
Then we get
\[
C_9+C_{10}\lesssim \varepsilon^2d_k^22^{-k}\left(||e^{\mathcal Rt'}(Q_{11},\varepsilon Q_{12})_\psi||_{\tilde L^2_t(B^\frac 12)}^4+||e^{\mathcal Rt'}\varepsilon\partial_x(Q_{11},\varepsilon Q_{12})_\psi||_{\tilde L^2_t(B^\frac 12)}^4+||e^{\mathcal Rt'}\partial_y(Q_{11},\varepsilon Q_{12})_\psi||_{\tilde L^2_t(B^\frac 12)}^4\right)
\]
so finally we have
\begin{equation}
\label{eqc}
\begin{aligned}
&\frac12||e^{\mathcal Rt'}(\varepsilon^2 \partial_x, \varepsilon\partial_y)\Delta^h_k(Q_{11},\varepsilon Q_{22})_\psi||^2_{L^\infty_t(L^2)}+a'\int_0^t{||e^{\mathcal Rt'}(\varepsilon^2 \partial_x, \varepsilon\partial_y)\Delta^h_k(Q_{11},\varepsilon Q_{22})_\psi(t')||^2_{L^2}dt'}\\
&+\lambda 2^k\int_0^t{\eta'(t)||e^{\mathcal Rt'}\Delta^h_k(\varepsilon^2 \partial_x, \varepsilon\partial_y)(Q_{11},\varepsilon Q_{22})_\psi(t')||_{L^2}^2 dt'}+\frac{\varepsilon^2}{2}\int_0^t{|| e^{\mathcal Rt'}\Delta^h_k(\Delta_\varepsilon Q_{11},\varepsilon \Delta_\varepsilon Q_{12})_\psi(t')||_{L^2}^2dt'}\\
&\le ||e^{a|D_x|}\Delta^h_k(\varepsilon^2\partial_x,\varepsilon\partial_y)((Q_{11})_0,\varepsilon (Q_{12})_0)||_{L^2}^2+ C\varepsilon^2d_k^22^{-k}||e^{\mathcal Rt'}\partial_x(Q_{11},\varepsilon Q_{12})_\psi||^2_{\tilde L^2_{t,\eta'(t)}(B^\frac 12)}\\
&+Cd_k^22^{-k}  ||e^{\mathcal Rt'}(Q_{11}, \varepsilon Q_{12})_\psi||^2_{\tilde L^2_{t,\eta'(t)}(B^\frac 12)}+C \varepsilon^4d_k^22^{-k}||e^{\mathcal Rt'}\partial_xu_\psi||^2_{\tilde L^2_{t,\eta'(t)}(B^\frac 12)}\\
&+C\varepsilon^2d_k^22^{-k}||e^{\mathcal Rt'}(Q_{11},\varepsilon Q_{12})_\psi||_{\tilde L^2_t(B^\frac 12)}^4+C\varepsilon^2d_k^22^{-k}||e^{\mathcal Rt'}\varepsilon\partial_x(Q_{11},\varepsilon Q_{12})_\psi||_{\tilde L^2_t(B^\frac 12)}^4\\
&+ C\varepsilon^4d_k^22^{-k}||e^{\mathcal Rt'}\partial_xu_\psi||^2_{\tilde L^2_{t,\eta'(t)}(B^\frac 12)}+C\varepsilon^2d_k^22^{-k}||e^{\mathcal Rt'}\partial_y(Q_{11},\varepsilon Q_{12})_\psi||_{\tilde L^2_t(B^\frac 12)}^4
\end{aligned}
\end{equation}

Now we come back to prove the main Theorem.
\begin{proof}
Multiply the inequalities above $2^k$, take square root and sum up over $\mathbb Z$, and recall that $a>0$. Then we have
\begin{equation}
\label{eqabc}
\begin{aligned}
&||e^{\mathcal Rt'}(\varepsilon u,\varepsilon^2v)_\psi||_{\tilde L^\infty_t(B^\frac 12)}+||e^{\mathcal Rt'}(Q_{11},\varepsilon Q_{12})_\psi||_{\tilde L^\infty_t(B^\frac 12)}+||e^{\mathcal Rt'}(\varepsilon^2\partial_x,\varepsilon\partial_y)(Q_{11},\varepsilon Q_{12})_\psi||_{\tilde L^\infty_t(B^\frac 12)}\\
&+\sqrt\lambda||e^{\mathcal Rt'}(\varepsilon u,\varepsilon^2 v)_\psi||_{\tilde L^2_{t,\eta'(t)}(B^1)}+c||e^{\mathcal Rt'}\partial_y(\varepsilon u,\varepsilon^2v)_\psi||_{\tilde L^2_t(B^\frac 12)}+c||e^{\mathcal Rt'}\partial_x(\varepsilon u,\varepsilon^2v)_\psi||_{\tilde L^2_t(B^\frac 12)}\\
&+\sqrt\lambda||e^{\mathcal Rt'}(Q_{11},\varepsilon Q_{12})_\psi||_{\tilde L^2_{t,\eta'(t)}(B^1)}+c||e^{\mathcal Rt'}\partial_y(Q_{11},\varepsilon Q_{12})_\psi||_{\tilde L^2_t(B^\frac 12)}+c||e^{\mathcal Rt'}\varepsilon\partial_x(Q_{11},\varepsilon Q_{12})_\psi||_{\tilde L^2_t(B^\frac 12)}\\
&+\sqrt\lambda||e^{\mathcal Rt'}(\varepsilon^2\partial_x,\varepsilon\partial_y)(Q_{11},\varepsilon Q_{12})_\psi||_{\tilde L^2_{t,\eta'(t)}(B^1)}+c\varepsilon||e^{\mathcal Rt'}\Delta_\varepsilon(Q_{11},\varepsilon Q_{12})_\psi||_{\tilde L^2_t(B^\frac 12)}\\
&\le||e^{a|D_x|}(\varepsilon u_0,\varepsilon^2v_0)||_{B^\frac 12}+||e^{a|D_x|}((Q_{11})_0,\varepsilon (Q_{12})_0)||_{B^\frac 12}+||e^{a|D_x|}(\varepsilon^2\partial_x,\varepsilon\partial_y)((Q_{11})_0,\varepsilon (Q_{12})_0)||_{B^\frac 12}\\
&+C\varepsilon^{\frac 32}||e^{\mathcal Rt'}(\varepsilon u,\varepsilon^2 v)_\psi||_{\tilde L^2_{t,\eta'(t)}(B^1)}+C\varepsilon^{\frac 12}||e^{\mathcal R't'}(\varepsilon^2\partial_x,\varepsilon\partial_y)(Q_{11},\varepsilon Q_{12})_\psi||_{\tilde L^2_{t,\eta'(t)}(B^1)}\\
&+C||e^{\mathcal Rt'}\partial_y(Q_{11},\varepsilon Q_{12})_\psi||_{\tilde L^2_{t,\eta'(t)}(B^\frac12)}+C||e^{\mathcal Rt'}\partial_x (\varepsilon u,\varepsilon^2 v)_\psi||_{\tilde L^2_{t,\eta'(t)}(B^\frac12)}+ C ||e^{\mathcal R't'} \varepsilon\partial_x(Q_{11}, \varepsilon Q_{12})_\psi||_{\tilde L^2_{t,\eta'(t)}(B^\frac12)}\\
&+C\varepsilon ||e^{\mathcal Rt'}\Delta_\varepsilon(Q_{11},\varepsilon Q_{12})_\psi||_{\tilde L^2_{t,\eta'(t)}(B^\frac12)}+c||e^{\mathcal R't'}\partial_y(Q_{11},\varepsilon Q_{12})_\psi||^2_{\tilde L^2_t(B^\frac 12)}+c||e^{\mathcal R't'}\partial_x(Q_{11},\varepsilon Q_{12})_\psi||^2_{\tilde L^2_t(B^\frac 12)}\\
\end{aligned}
\end{equation}
Recall that on $(0,T^*]$, we have $||\partial_y(Q_{11},\varepsilon Q_{12})_\psi||_{\tilde L^2_t(B^\frac 12)}+||\varepsilon\partial_x(Q_{11},\varepsilon Q_{12})_\psi||_{\tilde L^2_t(B^\frac 12)}<\eta'(t)<\delta$. Choosing $\lambda\gg C^2$ and $\delta$ small enough, we have
\begin{equation}
\label{eqabc1}
\begin{aligned}
&||e^{\mathcal Rt'}(\varepsilon u,\varepsilon^2v)_\psi||_{\tilde L^\infty_t(B^\frac 12)}+||e^{\mathcal Rt'}(Q_{11},\varepsilon Q_{12})_\psi||_{\tilde L^\infty_t(B^\frac 12)}+||e^{\mathcal Rt'}(\varepsilon^2\partial_x,\varepsilon\partial_y)(Q_{11},\varepsilon Q_{12})_\psi||_{\tilde L^\infty_t(B^\frac 12)}\\
&+\frac c2||e^{\mathcal Rt'}\partial_y(\varepsilon u,\varepsilon^2v)_\psi||_{\tilde L^2_t(B^\frac 12)}+\frac c2||e^{\mathcal Rt'}\partial_x(\varepsilon u,\varepsilon^2v)_\psi||_{\tilde L^2_t(B^\frac 12)}\\&+\frac c2||e^{\mathcal Rt'}\partial_y(Q_{11},\varepsilon Q_{12})_\psi||_{\tilde L^2_t(B^\frac 12)}+\frac c2||e^{\mathcal Rt'}\varepsilon\partial_x(Q_{11},\varepsilon Q_{12})_\psi||_{\tilde L^2_t(B^\frac 12)}+\frac c2\varepsilon||e^{\mathcal R't'}\Delta_\varepsilon(Q_{11},\varepsilon Q_{12})_\psi||_{\tilde L^2_t(B^\frac 12)}\\
&\le||e^{a|D_x|}(\varepsilon u_0,\varepsilon^2v_0)||_{B^\frac 12}+||e^{a|D_x|}((Q_{11})_0,\varepsilon (Q_{12})_0)||_{B^\frac 12}+||e^{a|D_x|}(\varepsilon^2\partial_x,\varepsilon\partial_y)((Q_{11})_0,\varepsilon (Q_{12})_0)||_{B^\frac 12}:=c_0\\
\end{aligned}
\end{equation}
so for $t\le T^*$, we have
\begin{equation*}
\begin{aligned}
||\partial_y(\varepsilon u,\varepsilon^2v)_\psi||_{\tilde L^2_t(B^\frac 12)}&+||\partial_x(\varepsilon u,\varepsilon^2v)_\psi||_{\tilde L^2_t(B^\frac 12)}\\&+||\partial_y(Q_{11},\varepsilon Q_{12})_\psi||_{\tilde L^2_t(B^\frac 12)}+||\varepsilon\partial_x(Q_{11},\varepsilon Q_{12})_\psi||_{\tilde L^2_t(B^\frac 12)}\le \frac{2c_0}{c}e^{-\mathcal Rt}\\
\end{aligned}
\end{equation*}
Recall the definition of $T^*$. If we choose $c_0$ small enough, such that $c_0\le \frac{c\delta}{4}$, then we deduce that $T^*$ equals to $\infty$. Thus we  have proved the Theorem.
\end{proof}
\section{Global well-posedness of the hydrostatic system and the convergence}\label{Sechydro}
In this section, we study the global well-posedness of the hydrostatic approximate equations \eqref{dim2limit1} with small analytic data and justify the limit from the scaled anisotropic system to the hydrostatic system.
\subsection{Global solutions of the hydrostatic approximate}
Recall that for the hydrostatic limit equation \eqref{q11l2}, we should skip the ordinary case that $\partial_yU=0$. Then \eqref{q11l2} becomes $Q_{11}=Q_{12}=0$. So we only need to focus on the equation \eqref{dim2limit1}. Recall that $U$ has the expression given in  \eqref{shearflowu} and we suppose that $\sum_{m\ge 0}m|c_*(m)|<\mathsf c^*$. Define 
\begin{equation}
\label{hydro1}
u_{\phi}(t,x,y)=\mathcal F^{-1}_{\xi\to x}(e^{\phi(t,\xi)}\hat u(t,\xi,y)) \quad\text{with}\quad \phi(t,\xi):=(a-\lambda \theta(t))|\xi|
\end{equation}
where $\theta(t)$ denotes the evolution of the analytic band of $u$, determined by
\begin{equation}
\label{hydro2}
\theta'(t)=\sum_{m>0}m|c_*(m)|e^{-m^2\pi^2t}\quad\mbox{with}\quad \theta(0)=0\end{equation}
Notice that $\theta(t)\le\sum_{m>0} \frac{|c_*(m)|}{m\pi^2}$ for any $t>0$. We first prove the following proposition.
\begin{prop}
\label{propdim2lim1}
For any $s>0$, there exists some constant $\mathsf c^*>0$, such that if $\sum_{m>0}\frac{|c_*(m)|}{m}<\mathsf c^*$, then
\begin{equation}
\label{dimension2lim5}
||e^{\mathcal Rt'}u_\phi||_{\tilde L^\infty_t(B^s)}+||e^{\mathcal Rt'}\partial_yu_\phi||_{\tilde L^2_t(B^s)}\le ||e^{a|D_x|}u_0||_{B^s}
\end{equation}
\end{prop}
\begin{proof}
Recall that $|D_x|$ denotes the Fourier multiplier with symbol $|\xi|$. Using \eqref{dim2limit1}, we deduce that $u_\phi$ satisfies 
\begin{equation}
\label{dimension2lim1}
\partial_tu_\phi+\lambda\theta'(t)|D_x|u_\psi+(U\partial_xu)_\phi+(v\partial_yU)
_\phi-\partial_y^2u_\phi+\partial_xp_\phi=0\end{equation}
Appling $\Delta^h_k$ to \eqref{dimension2lim1} and taking the $L^2$ inner product with $\Delta^h_ku_\phi$, we have
\begin{equation}
\label{dimension2lim2}
\begin{aligned}
&\frac 12\frac {d}{dt}||\Delta^h_ku_\phi||^2_{L^2}
+\lambda\theta'(|D_x|\Delta^h_ku_\phi,\Delta^h_ku_\phi)_{L^2}+||\Delta^h_k\partial_yu_\phi||^2_{L^2}\\
&=-(\Delta ^h_k(U\partial_xu)_\phi,\Delta^h_ku_\phi)_{L^2}-(\Delta ^h_k(v\partial_yU)_\phi,\Delta^h_ku_\phi)_{L^2}-(\Delta ^h_k(\partial_xp_\phi,\Delta^h_ku_\phi)_{L^2}
\end{aligned}
\end{equation}
Similarly as the anisotropic case, we obtain that $(\Delta ^h_k(U\partial_xu)_\phi,\Delta^h_ku_\phi)_{L^2}=0$, and
\[
(\Delta ^h_k(\partial_xp_\phi,\Delta^h_ku_\phi)_{L^2}=-(\Delta ^h_k(p_\phi,\Delta^h_k\partial_xu_\phi)_{L^2}=-(\Delta ^h_k(p_\phi,\Delta^h_k\partial_yv_\phi)_{L^2}=-(\Delta ^h_k(\partial_yp_\phi,\Delta^h_kv_\phi)_{L^2}=0.
\]
Next, recalling that $\partial_yU\lesssim \theta'(t)$,  we can directly deduce that for any $s>0$,
\begin{equation}
\label{dimension2lim3s}
-(\Delta ^h_k(v\partial_yU)_\phi,\Delta^h_ku_\phi)_{L^2}\lesssim d_k^22^{-2ks}||u_\phi||^2_{\tilde L^2_{t,\theta'(t)}(B^{s+\frac 12})}
\end{equation}
Multiplying  \eqref{dimension2lim2} by $e^{2\mathcal Rt}$ and integrating over $[0,t]$, we obtain that there exists a constant $C>0$ such that
\begin{equation}
\label{dimension2lim3}
\begin{aligned}
&\frac{1}{2}||e^{\mathcal Rt'}\Delta^h_ku_\phi||^2_{L^\infty_t(L^2)}+\lambda2^k\int_0^t{\theta'(t')||e^{\mathcal Rt'}\Delta^h_ku_\phi(t')||^2_{L^2}dt'}+\frac12 ||e^{\mathcal Rt'}\Delta^h_k\partial_yu_\phi||^2_{L^2_t(L^2)}\\&\le ||e^{a|D_x|}\Delta^h_ku_0||^2_{L^2}+C d_k^22^{-2ks}||e^{\mathcal Rt'}u_\phi||^2_{\tilde L^2_{t,\theta'(t)}(B^{s+\frac 12})}
\end{aligned}
\end{equation}
For any $s>0$, we multiply \eqref{dimension2lim3} by $2^{2ks}$, take square root and sum up over $\mathbb Z$ to obtain that for any $t\le T^*$,
\begin{equation}
\label{dimension2lim4}
||e^{\mathcal Rt'}u_\phi||_{\tilde L^\infty_t(B^s)}+\sqrt\lambda ||e^{\mathcal Rt'}u_\phi||_{\tilde L^2_{t,\theta'(t)}(B^{s+\frac 12})}+||e^{\mathcal Rt'}\partial_yu_\phi||_{\tilde L^2_t(B^s)}\le ||e^{a|D_x|}u_0||_{B^s}+C||e^{\mathcal Rt'}u_\phi||_{\tilde L^2_{t,\theta'(t)}(B^{s+\frac 12})}
\end{equation}
Finally, if $\mathsf c^*<\frac{a\pi^2}{C^2}$, then we can choose a suitable $\lambda$ such that  $\lambda> C^2$, and we have proved the proposition.
\end{proof}
\begin{rmk}
Here we just need $\theta(t)$ to be small enough instead of $\theta'(t)$, which allows the condition on $U$ to be weaker.
\end{rmk}
Next, we prove
\begin{prop}
\label{propdim2lim2}
For any $s>0$,  there exists some constant $\mathsf c^*>0$, such that if $\sum_{m>0}m{|c_*(m)|}<\mathsf c^*$, then there exists a constant $C>0$ such that
\begin{equation}
\label{dimension2lim8}
||e^{\mathcal Rt'}\partial_yu_\phi||_{\tilde L^\infty_t(B^s)}+||e^{\mathcal Rt'}\partial^2_yu_\phi||_{\tilde L^2_t(B^s)}\le C(||e^{a|D_x|}u_0||_{B^s}+||e^{a|D_x|}u_0||_{B^{s+1}}+||e^{a|D_x|}\partial_yu_0||_{B^s})
\end{equation}
\end{prop}
\begin{proof}
We come back to the equation \eqref{dim2limit1}, and we follow the idea of Lemma 3.2 of \cite{PZZ}. Applying $\partial_y$ to both sides, we have 
\[
\partial_t\partial_yu+U\partial_x\partial_yu+v\partial_y^2U-\partial_y^3u+\partial_x\partial_yp=0
\]
Recall that
$(\Delta^h_k(U\partial_x\partial_yu)_\phi, \Delta^h_k\partial_yu_\phi)_{L^2}=0.$ Similarly as before, we obtain that for any $0<t<T^*$,
\begin{equation}
\label{dimension2lim4}
\begin{aligned}
&\frac 12||e^{\mathcal Rt'}\Delta^h_k\partial_yu_\phi||^2_{L^\infty_t(L^2)}+\lambda 2^k\int_0^t{\theta'(t)||e^{\mathcal Rt'}\Delta^h_k\partial_yu_\phi(t')||^2_{L^2}dt'}+\frac 12||e^{\mathcal Rt'}\Delta^h_k\partial^2_yu_\phi||^2_{L^2_t(L^2)}\\
&\le\frac12||e^{a|D_x|}\Delta^h_k\partial_yu_0||^2_{L^2}
+\int_0^t{|(e^{\mathcal Rt'}\Delta^h_k(v\partial_y^2U)_\psi, e^{\mathcal Rt'}\Delta^h_k\partial_yu_\phi|)_{L^2}|dt'}+\int_0^t{|(e^{\mathcal Rt'}\Delta^h_k\partial_xp_\phi, e^{\mathcal Rt'}\Delta^h_k\partial^2_yu_\phi|)_{L^2}|dt'}
\end{aligned}
\end{equation}
From integration by parts, we have 
\[(\Delta^h_k(v\partial_y^2U)_\phi, \Delta^h_k\partial_yu_\phi)_{L^2}=-(\Delta^h_k(v\partial_yU)_\phi, \Delta^h_k\partial_y^2u_\phi)_{L^2}+(\Delta^h_k(\partial_xu\partial_yU)_\phi, \Delta^h_k\partial_yu_\phi)_{L^2}
\]
Notice that $\partial_yU\lesssim \theta'(t)\le \mathsf c^*$, so
\begin{equation}
\label{extra1}
\begin{aligned}
(\Delta^h_k(v\partial_yU)_\phi, \Delta^h_k\partial_y^2u_\phi)_{L^2}&\lesssim  d_k^22^{-2ks}||\partial_y^2u_\phi||^2_{\tilde L^2_{t,\theta'(t)}(B^s)}+d_k^22^{-2ks}||u_\phi||^2_{\tilde L^2_{t,\theta'(t)}(B^{s+1})}\\
&\lesssim  \mathsf c^*d_k^22^{-2ks}||\partial_y^2u_\phi||^2_{\tilde L^2_{t}(B^s)}+\mathsf c^*d_k^22^{-2ks}||u_\phi||^2_{\tilde L^2_{t}(B^{s+1})}\\
\end{aligned}
\end{equation}
and similarly from Poincar\'e inequality,
\begin{equation}
\label{extra2}
\begin{aligned}
(\Delta^h_k(\partial_xu\partial_yU)_\phi, \Delta^h_k\partial_yu_\phi)_{L^2}
&\lesssim  d_k^22^{-2ks}||\partial_yu_\phi||_{\tilde L^2_{t,\theta'(t)}(B^{s+\frac 12})}||u_\phi||_{\tilde L^2_{t, \theta'(t)}(B^{s+\frac 12})}
\lesssim  \frac{1}{\mathcal R}d_k^22^{-2ks}||\partial_yu_\phi||^2_{\tilde L^2_{t,\theta'(t)}(B^{s+\frac 12})}\\
\end{aligned}
\end{equation}
Next, for the $\partial_xp_\phi$ term, we integrate \eqref{dim2limit1} for $y\in [0,1]$ to get that

\begin{equation}
\label{dimension2lim5}
\partial_xp=\partial_yu(t,x,1)-\partial_yu(t,x,0)-2\partial_x\int_0^1{U(t,y)u(t,x,y)dy}
\end{equation}
so we have
\begin{equation*}
\begin{aligned}
(e^{\mathcal Rt}\Delta^h_k\partial_xp_\phi(t), e^{\mathcal Rt}\Delta^h_k\partial^2_yu_\phi(t))_{L^2}
&=\int_{\mathbb R}{e^{\mathcal Rt}\Delta^h_k\partial_xp_\phi\cdot e^{\mathcal Rt}(\Delta^h_k\partial_yu_\phi(t,x,1)-\Delta^h_k\partial_yu_\phi(t,x,0))dx}\\
&=\int_{\mathbb R}{ (e^{\mathcal Rt}(\Delta^h_k\partial_yu_\phi(t,x,1)-\Delta^h_k\partial_yu_\phi(t,x,0)))^2dx}
\\
&+2\int_{\mathbb R}{\left(e^{\mathcal Rt}\int_0^1{\Delta^h_k\partial_x(Uu)_\phi dy}\right)\cdot e^{\mathcal Rt}(\Delta^h_k\partial_yu_\phi(t,x,1)-\Delta^h_k\partial_yu_\phi(t,x,0))dx}\\
&\lesssim ||e^{\mathcal Rt}\Delta^h_k\partial_yu_\phi(t)||^2_{L^\infty_v(L^2_h)}+||e^{\mathcal Rt}\Delta^h_k\partial_x(Uu)_\phi(t)||^2_{L^1_v(L^2_h)}
\end{aligned}
\end{equation*}
Notice that for any $s>0$,
\begin{equation*}
\begin{aligned}
||e^{\mathcal Rt}\Delta^h_k\partial_x(Uu)_\phi(t)||_{L^2_t({L^1_v(L^2_h))}}&\le ||e^{\mathcal Rt}\Delta^h_k\partial_x(Uu)_\phi(t)||_{L^2_t(L^2)}\\&\lesssim d_k2^{-ks}||U||_{L^2_tL^\infty_y}||e^{\mathcal Rt'}u_\phi||_{\tilde L^\infty_{t}(B^{s+1})} \lesssim \mathsf c^*d_k2^{-ks}||e^{\mathcal Rt'}u_\phi||_{\tilde L^\infty_{t}(B^{s+1})}
\end{aligned}
\end{equation*}
and similarly as the proof of Proposition 4.2 of \cite{PZZ}, we have
\[
||e^{\mathcal Rt}\Delta^h_k\partial_yu_\phi||^2_{L^\infty_v(L^2_h)}\le 2||e^{\mathcal Rt}\Delta^h_k\partial_yu_\phi||_{L^2}||e^{\mathcal Rt}\Delta^h_k\partial^2_yu_\phi||_{L^2}
\]
Combining the estimates together, we have
\begin{equation}
\label{extra3}
\begin{aligned}
&\int_0^t{|(e^{\mathcal Rt}\Delta^h_k\partial_xp_\phi(t), e^{\mathcal Rt}\Delta^h_k\partial^2_yu_\phi(t))_{L^2}|dt'}\\&\le \frac 14||e^{\mathcal Rt}\Delta^h_k\partial^2_yu_\phi||_{L^2_t(L^2)}^2+ C||e^{\mathcal Rt}\Delta^h_k\partial_yu_\phi||_{L^2_t(L^2)}^2+C(\mathsf c^*)^2d^2_k2^{-2ks}||e^{\mathcal Rt'}u_\phi||_{\tilde L^\infty_{t}(B^{s+1})}^2\\
&\le\frac 14||e^{\mathcal Rt}\Delta^h_k\partial^2_yu_\phi||_{L^2_t(L^2)}^2+Cd^2_k2^{-2ks}(||e^{\mathcal Rt}\Delta^h_k\partial_yu_\phi||_{\tilde L^2_t(B^s)}^2+(\mathsf c^*)^2||e^{\mathcal Rt'}u_\phi||_{\tilde L^\infty_{t}(B^{s+1})}^2)\\
\end{aligned}
\end{equation}
Take the square root and sum up in $\mathbb Z$, so from \eqref{extra1} to \eqref{extra3}, we have
\begin{equation*}
\begin{aligned}
&||e^{\mathcal Rt'}\partial_yu_\phi||_{\tilde L^\infty_t(B^s)}+\sqrt\lambda||e^{\mathcal Rt'}\partial_yu_\phi||_{\tilde L^2_{t,\theta'(t)}(B^{s+\frac 12})}+||e^{\mathcal Rt'}\partial^2_yu_\phi||_{\tilde L^2_t(B^s)}\\
&\le||e^{a|D_x|}\partial_yu_0||_{B^s}+C(||e^{\mathcal Rt'}\partial_yu_\phi||_{\tilde L^2_t(B^s)}+\mathsf c^*||e^{\mathcal Rt'}\partial_yu_\phi||_{\tilde L^\infty_{t}(B^{s+1})}+\frac {1}{ \mathcal R}||e^{\mathcal Rt'}\partial_yu_\phi||^2_{\tilde L^2_{t,\theta'(t)}(B^{s+\frac 12})}
)\\
&+C\mathsf c^*(||e^{\mathcal Rt}\partial^2_yu_\phi||_{\tilde L^2_t(B^s)}+||e^{\mathcal Rt'}\partial_yu_\phi||_{\tilde L^2_t(B^{s+1})})
\end{aligned}
\end{equation*}
together with the Proposition \ref{propdim2lim1} and Poincare inequality, if $\mathsf c^*$ is small enough, we choose $\lambda\gg\frac{C^2}{\mathcal R^2}$ and the proposition is proved.
\end{proof}
Now we come back to finish proving Theorem 1.2.
\begin{proof}
The first two inequalities have been proved in the two previous propositions. The main structure of proving the third one is similar to Theorem 1.2 of \cite{PZZ}. For the equation \eqref{dim2limit1}, we apply the operator $\Delta^h_k$ and take the $L^2$ inner product with $e^{2\mathcal Rt}\Delta^h_k(\partial_tu)_\phi$ to obtain that
\begin{equation}
\label{thm12extra0}
\begin{aligned}
||e^{\mathcal Rt}\Delta^h_k(\partial_tu)_\phi||^2_{L^2}&=e^{2\mathcal Rt}(\Delta^h_k\partial^2_yu_\phi, \Delta^h_k(\partial_yu)_\phi)_{L^2}\\&-e^{2\mathcal Rt}(\Delta^h_k(U\partial_xu)_\phi, \Delta^h_k(\partial_yu)_\phi)_{L^2}-e^{2\mathcal Rt}(\Delta^h_k(v\partial_yU)_\phi, \Delta^h_k(\partial_yu)_\phi)_{L^2}
\end{aligned}
\end{equation}
From integration by parts we get
\begin{equation}
\label{thm12extra1}
(e^{2\mathcal Rt}\Delta^h_k\partial^2_yu_\psi, \Delta^h_k(\partial_yu)_\psi)_{L^2}=-\frac 12\frac {d}{dt}||e^{\mathcal Rt}\Delta^h_h\partial_yu_\psi||^2_{L^2}-\theta'(t)2^k||e^{\mathcal Rt}\Delta^h_k\partial_yu_\psi||^2_{L^2}
\end{equation}
For an arbitary $s>0$, we multiply \eqref{thm12extra0} by $2^{2ks}$, use \eqref{thm12extra1} and sum over $\mathbb Z$, to get
\begin{equation*}
\begin{aligned}
||e^{\mathcal Rt'}(\partial_tu)_\psi||_{\tilde L^2_t(B^s)}&+||e^{\mathcal Rt'}\partial_yu_\psi||_{\tilde L^\infty_t(B^s)}\\&\le C(||e^{a|D_x|}\partial_yu_0||_{B^s}+||e^{\mathcal Rt'}(U\partial_xu)_\psi||_{\tilde L^2_t(B^s)}+||e^{\mathcal Rt'}(v\partial_yU)_\psi||_{\tilde L^2_t(B^s)})
\end{aligned}
\end{equation*}
Recalling that $U$ and $\partial_yU$ are uniformly bounded by $\mathsf c^*$, by Poincare inequality, we have
\[
||e^{\mathcal Rt'}(U\partial_xu)_\phi||_{\tilde L^2_t(B^s)}, ||e^{\mathcal Rt'}(v\partial_yU)_\phi||_{\tilde L^2_t(B^s)}\lesssim ||e^{\mathcal Rt'}\partial_yu_\phi||_{\tilde L^2_t(B^{s+1})}
\]
so finally we obtain
\[
||e^{\mathcal Rt'}(\partial_tu)_\phi||_{\tilde L^2_t(B^s)}+||e^{\mathcal Rt'}\partial_yu_\phi||_{\tilde L^\infty_t(B^s)}\le C(||e^{a|D_x|}\partial_yu_0||_{B^s}+||e^{a|D_x|}u_0||_{\tilde L^2_t(B^{s+1})})
\]
and the theorem is proved.
\end{proof}
\begin{rmk}
Notice that in this subsection, the evolution equation $\theta(t)$ depends only on $U$. In fact, in next subsection, we need to add the another term $||\partial_yu_\phi(t')||_{B^\frac 12}$ to prove the convergence. Notice that
\[
\int_0^t||\partial_yu_\phi(t')||_{B^\frac 12}dt'\lesssim ||e^{\mathcal Rt'}\partial_yu_\phi||_{\tilde L^2_t(B^\frac 12)}\lesssim ||e^{a|D_x|}u_0||_{B^\frac 12}
\]
so if we also suppose $||e^{a|D_x|}u_0||_{B^\frac 12}$ is small enough, then $||\partial_yu_\phi(t')||_{B^\frac 12}$ will also appear in $\theta(t)$. Comparing with the anistropic case, we can similarly define the evolution of the analytic band of $u^\varepsilon-u$ to prove the convergence. The details can be seen later.
\end{rmk}

\subsection{Convergence to the hydrostatic system} Define $w^1_\varepsilon:=u^\varepsilon-u, w^2_\varepsilon:=v^\varepsilon-v$. Then they satisfy the equations:

\begin{equation}
\label{dim2limu1}
\begin{cases}
\varepsilon\partial_tw^1_{\varepsilon}-\varepsilon^3\partial_x^2w^1_\varepsilon-\varepsilon\partial_y^2w^1_\varepsilon+\varepsilon\partial_xq_\varepsilon+\partial_xR_{11}+\partial_yR_{21}\\
=\varepsilon^3\partial_x^2u_\varepsilon-\varepsilon[(U+ \varepsilon u^{\varepsilon})\partial_xu^{\varepsilon}-U\partial_xu]-\varepsilon [v^{\varepsilon}\partial_y(U+\varepsilon  u^{\varepsilon})-v\partial_yU]\\
\varepsilon^2\partial_tw^2_{\varepsilon}-\varepsilon^4\partial_x^2w^2_\varepsilon-\varepsilon^2 \partial_y^2w^2_\varepsilon+\partial_yq^{\varepsilon}+\partial_xR_{12}+\partial_yR_{22}\\=-\varepsilon^2\partial_tv+\varepsilon^4\partial_x^2v+\varepsilon^2 \partial_y^2v -\varepsilon^2(U+\varepsilon u^{\varepsilon})\partial_xv^{\varepsilon}-\varepsilon^3 v^{\varepsilon}\partial_yv^{\varepsilon}\\
\partial_xw^1_{\varepsilon}+\partial_yw^2_{\varepsilon}=0
\end{cases}
\end{equation}

 For a function $u$, we define 
\begin{equation}
\label{functiontheta}
u_\Theta(t,x,y):=\mathcal F^{-1}_{\xi\to x}(e^{\Theta(t,\xi)}\hat u(t,\xi,y)), \quad \Theta(t,\xi):=(a-\mu\zeta(t))|\xi|
\end{equation}
where $\mu\ge\lambda$ will be determined later, and $\zeta(t)$ is given by
\[
\zeta'(t)=||(\varepsilon\partial_x,\partial_y)u^\varepsilon_\psi(t')||_{B^\frac 12}+||(\varepsilon\partial_x,\partial_y)(Q_{11}^\varepsilon, \varepsilon Q_{12}^\varepsilon)_\psi(t')||_{B^\frac 12}+\sum_{m>0}m|c_m|e^{-m^2\pi^2t}+||u_{\phi}(t')||_{B^\frac 12}, \quad \zeta(0)=0
\]

Out of Theorems \ref{thm11} and \ref{thm12}, we deduce that
\[
||\varepsilon u_\Theta||_{\tilde L^\infty(\mathbb R^+, B^\frac 12)}+||u_\Theta||_{\tilde L^\infty(\mathbb R^+, B^\frac 12\bigcap B^\frac 52)}+||\partial_y u_\Theta||_{\tilde L^2(\mathbb R^+, B^\frac 12\bigcap B^\frac 52)}+||(\partial_tu)_\Theta||_{\tilde L^2(\mathbb R^+, B^\frac 32)}\le M
\]
where $M>1$ is a constant independant of $\varepsilon$. For convenience, we define 
\[
\mathcal S^1:=\varepsilon^3\partial_x^2u_\varepsilon-\varepsilon[(U+ \varepsilon u^{\varepsilon})\partial_xu^{\varepsilon}- U\partial_xu]-\varepsilon [v^{\varepsilon}\partial_y(U+\varepsilon  u^{\varepsilon})-v\partial_yU]
\]
and
\[
\mathcal S^2:=-\varepsilon^2\partial_tv+\varepsilon^4\partial_x^2v+\varepsilon^2 \partial_y^2v -\varepsilon^2(U+\varepsilon u^{\varepsilon})\partial_xv^{\varepsilon}-\varepsilon^3 v^{\varepsilon}\partial_yv^{\varepsilon}
\]
We prove two propositions about $\mathcal S^1$ and $\mathcal S^2$:
\begin{prop}
\label{proplim1}
\begin{equation}
\label{dim2limeq3}
\begin{aligned}
&\int_0^t{|(\Delta^h_k\mathcal S^1_\Theta,\Delta^h_kw^1_\Theta)_{L^2}|dt'}\lesssim \varepsilon d_k^22^{-k}||(\varepsilon w^1)_\Theta||_{\tilde L^2_{t, \zeta'(t)}(B^1)}||u_\Theta||^\frac 12_{\tilde L^\infty_t(B^\frac 32)}||\partial_y(\varepsilon w^1)_\Theta||_{\tilde L^2_{t}(B^{\frac 12})}\\
&+\varepsilon^3 d_k^22^{-k}||\partial_yu_\Theta||_{\tilde L^2_t(B^\frac 32)}||(\varepsilon w^1)_\Theta||_{\tilde L^2_t(B^\frac 32)}+ d_k^22^{-k}||(\varepsilon w^1)_\Theta||^2_{\tilde L^2_{t, \zeta'(t)}(B^1)}\\
&+\varepsilon d_k^22^{-k}(||u_\Theta||_{\tilde L^\infty_t(B^\frac 12) } ||\partial_yu_\Theta||_{\tilde L^2_t(B^{\frac 32}) }+||u_\Theta||_{\tilde L^\infty_t(B^{\frac 32}) } ||\partial_yu_\Theta||_{\tilde L^2_t(B^\frac 12) })||(\varepsilon w^1)_\Theta||_{\tilde L^2_{t}(B^{\frac 12})}\\
\end{aligned}
\end{equation}
\end{prop}
\begin{proof}
From integration by parts and Poincare inequality, we have
\begin{equation}
\label{dim2limeq4}
\varepsilon^4\int_0^t{|(\Delta^h_k\partial_x^2u_\Theta, \Delta^h_k w^1_\Theta)_{L^2}|dt'}\lesssim \varepsilon^3 d_k^22^{-k}||\partial_yu_\Theta||_{\tilde L^2_t(B^\frac 32)}||(\varepsilon w^1)_\Theta||_{\tilde L^2_t(B^\frac 32)}
\end{equation}
We next write 
\[
(U+\varepsilon u^\varepsilon)\partial_xu^\varepsilon-U\partial_xu=U\partial_xw^1 +u^\varepsilon\partial_xw^1+\varepsilon w^1\partial_xu+\varepsilon u\partial_xu\]
and
\[
v^\varepsilon\partial_y(U+\varepsilon u^\varepsilon)-v\partial_yU=w^2\partial_yU+\varepsilon w^2\partial_yu+\varepsilon v^\varepsilon\partial_yw^1+\varepsilon v\partial_yu
\]
Similarly as in Lemma 3.1 of \cite{PZZ} and recalling that  $U, \partial_yU\lesssim \zeta'(t)$, we have
\begin{equation}
\label{dimlimeq7}
\varepsilon^2\int_0^t{|(\Delta^h_k[(U+\varepsilon u^\varepsilon)\partial_xw^1+\partial_yUw^2]_\Theta, \Delta^h_k w^1_\Theta)_{L^2}|dt'}\lesssim d_k^22^{-k}||(\varepsilon w^1)_\Theta||^2_{\tilde L^2_{t, \zeta'(t)}(B^1)}
\end{equation}
Analogously  to (5.11) of \cite{PZZ} we get
\begin{equation}
\label{dimlimeq8}
\begin{aligned}
&\varepsilon\int_0^t{|(\Delta^h_k(\varepsilon w^1\partial_xu)_\Theta, \Delta^h_k(\varepsilon w^1)_\Theta)_{L^2}|dt'}\\&\lesssim\varepsilon d_k^22^{-k}(||(\varepsilon w^1)_\Theta||^2_{\tilde L^2_{t, \zeta'(t)}(B^1)}+||(\varepsilon w^1)_\Theta||_{\tilde L^2_{t, \zeta'(t)}(B^1)}||u_\Theta||^\frac 12_{\tilde L^\infty_t(B^\frac 32)}||\partial_y(\varepsilon w^1)_\Theta||_{\tilde L^2_{t}(B^{\frac 12})})
\end{aligned}
\end{equation}
Recall that $v^\varepsilon=w^2+v$. Similarly as in (5.13) of \cite{PZZ}, we obtain
\begin{equation}
\label{dimlimeq9}
\begin{aligned}
&\varepsilon\int_0^t{|(\Delta^h_k[v^\varepsilon\partial_y(\varepsilon w^1)]_\Theta, \Delta^h_k(\varepsilon w^1)_\Theta)_{L^2}|dt'}
\lesssim\varepsilon d_k^22^{-k}||(\varepsilon w^1)_\Theta||_{\tilde L^2_{t, \zeta'(t)}(B^1)}||u_\Theta||^\frac 12_{\tilde L^\infty_t(B^\frac 32)}||\partial_y(\varepsilon w^1)_\Theta||_{\tilde L^2_{t}(B^{\frac 12})}
\end{aligned}
\end{equation}
Finally, similarly as  (5.14) of \cite{PZZ},, we have
\begin{equation}
\label{dimlimeq10}
\begin{aligned}
&\varepsilon\int_0^t{|(\Delta^h_k(\varepsilon w^2\partial_yu)_\Theta, \Delta^h_k(\varepsilon w^1)_\Theta)_{L^2}|dt'}\lesssim\varepsilon d_k^22^{-k}||(\varepsilon w^1)_\Theta||^2_{\tilde L^2_{t, \zeta'(t)}(B^1)}\\&
\end{aligned}
\end{equation}
A  difference from \cite{PZZ} is that in here  there is an extra term $\varepsilon (u\partial_xu+v\partial_yu)$. Notice that for any $s>0$, from the law of product in anisotropic Besov space (see \cite{CPZ} for the proof) and Poincare inequality,
we have
\begin{equation}
\label{dimlimeq6}
||(u\partial_xu)_\Theta||_{\tilde L^2_t(B^s) }+||(v\partial_yu)_\Theta||_{\tilde L^2_t(B^s) }\lesssim ||u_\Theta||_{\tilde L^\infty_t(B^\frac 12) } ||\partial_yu_\Theta||_{\tilde L^2_t(B^{s+1}) }+||u_\Theta||_{\tilde L^\infty_t(B^{s+1}) } ||\partial_yu_\Theta||_{\tilde L^2_t(B^\frac 12) }
\end{equation}

and the proposition is proved after combining \eqref{dim2limeq4}-\eqref{dimlimeq6}.
\end{proof}

\begin{prop}
\label{proplim2}
\begin{equation}
\label{dim2limeq3extra}
\begin{aligned}
&\int_0^t{|(\Delta^h_k\mathcal S^2_\Theta,\Delta^h_kw^2_\Theta)_{L^2}|dt'}\lesssim \varepsilon d_k^22^{-k}||(\varepsilon w^1, \varepsilon^2 w^2)_\Theta||^2_{\tilde L^2_{t,\zeta'(t)}(B^1)}+\varepsilon^2 d_k^22^{-k}||\varepsilon^2w^2_\Theta||_{\tilde L^2_{t,\zeta'(t)}(B^1)}||\partial_yu_\Theta||_{\tilde L^2_{t}(B^{2})}\\&+ \varepsilon^2d_k^22^{-k}\left[\left(||(\partial_tu)_\Theta||_{\tilde L^2_t(B^\frac 32)}+||
\partial_yu_\Theta||_{\tilde L^2_t(B^\frac 32)}\right)||\partial_yw^2_\Theta||_{\tilde L^2_t(B^\frac 12)}+\varepsilon^2||
\partial_yu_\Theta||_{\tilde L^2_t(B^\frac 52)}||w^2_\Theta||_{\tilde L^2_t(B^\frac 32)}\right]\\
&+ \varepsilon d_k^22^{-k}||\varepsilon^2w^2_\Theta||_{\tilde L^2_{t,\zeta'(t)}(B^1)}\left\{||u^\varepsilon_\Theta||^\frac 12_{L^\infty_t(B^\frac 12)}||\partial_yu_\Theta||_{\tilde L^2_t(B^2)}+||u_\Theta||^\frac 12_{L^\infty_t(B^\frac 32)}\left(||\partial_yu_\Theta||_{\tilde L^2_t(B^\frac 32)}+||\partial_y(\varepsilon^2w^2)_\Theta||_{\tilde L^2_t(B^\frac 12)}\right)\right\}
\end{aligned}
\end{equation}
\end{prop}
\begin{proof}
Similarly as (5.15) of \cite{PZZ},we deduce from Poincar\'e inequality that
\begin{equation}
\label{eqs201}
\begin{aligned}
&\int_0^t{|\Delta^h_k(-\varepsilon^2\partial_tv+\varepsilon^4\partial_x^2v+\varepsilon^2\partial^2_yv)_\Theta, \Delta^h_kw^2_\Theta|_{L^2}dt}\\&\lesssim \varepsilon^2d_k^22^{-k}\left[\left(||(\partial_tu)_\Theta||_{\tilde L^2_t(B^\frac 32)}+||
\partial_yu_\Theta||_{\tilde L^2_t(B^\frac 32)}\right)||\partial_yw^2_\Theta||_{\tilde L^2_t(B^\frac 12)}+\varepsilon^2||
\partial_yu_\Theta||_{\tilde L^2_t(B^\frac 52)}||w^2_\Theta||_{\tilde L^2_t(B^\frac 32)}\right]
\end{aligned}
\end{equation}
We decompose $v^\varepsilon=v+w^2$. Recall that $U\lesssim \zeta'(t)$, so we have
\begin{equation*}
\int_0^t{|\Delta^h_k(U\partial_xw^2)_\Theta, \Delta^h_k w^2_\Theta|_{L^2}dt'}=0, \varepsilon^2\int_0^t{|\Delta^h_k(U\partial_xv)_\Theta, \Delta^h_k \varepsilon^2w^2_\Theta|_{L^2}dt'}\lesssim\varepsilon^2 d_k^22^{-k}||\varepsilon^2w^2_\Theta||_{\tilde L^2_{t,\zeta'(t)}(B^1)}||\partial_yu_\Theta||_{\tilde L^2_{t}(B^{2})}
\end{equation*}
Next, similarly as (5.16) and (5.17) of \cite{PZZ}, we have
\begin{equation}
\label{eqs202}
\varepsilon^3\int_0^t{|\Delta^h_k(u^\varepsilon\partial_xv^\varepsilon)_\Theta, \Delta^h_k \varepsilon^2w^2_\Theta|_{L^2}dt'}\lesssim \varepsilon d_k^22^{-k}||\varepsilon^2w^2_\Theta||_{\tilde L^2_{t,\zeta'(t)}(B^1)}(||\varepsilon^2w^2_\Theta||_{\tilde L^2_{t,\zeta'(t)}(B^1)}+||u^\varepsilon_\Theta||^\frac 12_{L^\infty_t(B^\frac 12)}||\partial_yu_\Theta||_{\tilde L^2_t(B^2)})
\end{equation}
Note that $v^\varepsilon\partial_yv^\varepsilon=v\partial_yw^2+w^2\partial_yw^2+v\partial_yv+w^2\partial_yv$. Recall from Lemma 3.3 of \cite{PZZ} that
\[
\varepsilon^3\int_0^t{|(\Delta^h_k(w^2\partial_yw^2)_\Theta, \Delta^h_k\varepsilon^2w^2_\Theta)_{L^2}|dt'}\lesssim \varepsilon d_k^2||(\varepsilon w^1_\Theta, \varepsilon^2 w^2_\Theta)||^2_{\tilde L^2_{t,\zeta'(t)}(B^1)}
\]
Next, recall that $\partial_yv=-\partial_xu$. Similarly as \eqref{dimlimeq8}, we have
\begin{equation}
\label{dimlimeq12}
\begin{aligned}
&\varepsilon^3\int_0^t{|(\Delta^h_k(w^2\partial_xu)_\Theta, \Delta^h_k\varepsilon^2w^2_\Theta)_{L^2}|dt'}\lesssim\varepsilon d_k^22^{-k}||\varepsilon^2w^2_\Theta||_{\tilde L^2_{t,\zeta'(t)}(B^1)}(||\varepsilon^2w^2_\Theta||_{\tilde L^2_{t,\zeta'(t)}(B^1)}+||u_\Theta||^\frac 12_{L^\infty_t(B^\frac 32)}||\partial_y(\varepsilon^2w^2)_\Theta||_{\tilde L^2_t(B^\frac 12)})
\end{aligned}
\end{equation}
and
\begin{equation}
\label{dimlimeq13}
\begin{aligned}
&\varepsilon^3\int_0^t{|(\Delta^h_k(v\partial_yw^2)_\Theta, \varepsilon^2\Delta^h_kw^2_\Theta)_{L^2}|dt'}\lesssim \varepsilon d_k^22^{-k}||(\varepsilon^2 w^2)_\Theta||_{\tilde L^2_{t, \zeta'(t)}(B^1)}||u_\Theta||^\frac 12_{\tilde L^\infty_t(B^\frac 32)}||\partial_y(\varepsilon^2 w^2)_\Theta||_{\tilde L^2_{t}(B^{\frac 12})}
\end{aligned}
\end{equation}
Finally, again using the law product of anistropic Besov norms, we have
\begin{equation}
\label{dimlimeq14}
\begin{aligned}
\varepsilon^3\int_0^t{|(\Delta^h_k(v\partial_yv)_\Theta, \Delta^h_k\varepsilon^2w^2_\Theta)_{L^2}|dt'}\lesssim \varepsilon^3 d_k^22^{-k}||(\varepsilon^2 w^2)_\Theta||_{\tilde L^2_{t, \zeta'(t)}(B^1)}||u_\Theta||^\frac 12_{\tilde L^\infty_t(B^\frac 32)}||\partial_yu_\Theta||_{\tilde L^2_{t}(B^{\frac 32})}
\end{aligned}
\end{equation}
and the propostion is proved after combining \eqref{eqs201} to \eqref{dimlimeq14}.
\end{proof}

Finally, we come back to proving Theorem \ref{thm3}.
\begin{proof}
Recall that in the hydrostatic limit case we assume that $Q_{11}=Q_{12}=0$ and we still drop '$\varepsilon$' in $Q^\varepsilon_{11}$ and $Q^\varepsilon_{12}$, for simplicity, as the anisotropic case.  By using a similar derivation of the scaled anisotropic system and recalling that the terms $R_{ij}$ can be estimated similarly as in the Propositions \ref{propa5}-\ref{propa8},  we obtain that 
\begin{equation}
\label{thm3eq1}
\begin{aligned}
&||\Delta^h_k(\varepsilon w^1_\Theta, \varepsilon^2 w^2_\Theta)||^2_{L^\infty_t(L^2)}+\mu 2^k\int_0^t{\zeta'(t') ||\Delta^h_k((\varepsilon w^1)_\Theta, (\varepsilon^2 w^2)_\Theta))||^2_{L^2}dt'}+||\Delta^h_k(Q_{11}, \varepsilon Q_{12})_\Theta)||^2_{L^\infty_t(L^2)}\\
&+\int_0^t{ ||\Delta^h_k\partial_y((\varepsilon w^1)_\Theta, (\varepsilon^2 w^2)_\Theta))||^2_{L^2}+\varepsilon^22^{2k}||\Delta^h_k(\varepsilon w^1)_\Theta, (\varepsilon^2 w^2)_\Theta))||^2_{L^2}dt'}+\int_0^t{ ||\Delta^h_k\Delta_\varepsilon((Q_{11},\varepsilon Q_{12})_\Theta))||^2_{L^2}dt'}\\
&+\mu 2^k\int_0^t{\zeta'(t') ||\Delta^h_k(Q_{11},\varepsilon Q_{12})_\Theta)||^2_{L^2}dt'}+\int_0^t{ ||\Delta^h_k\partial_y((Q_{11},\varepsilon Q_{12})_\Theta))||^2_{L^2}+\varepsilon^22^{2k}||\Delta^h_k((Q_{11}, \varepsilon Q_{12})_\Theta))||^2_{L^2}dt'}\\
&+||\Delta^h_k(\varepsilon^2\partial_x,\varepsilon\partial_y)(Q_{11}, \varepsilon Q_{12})_\Theta)||^2_{L^\infty_t(L^2)}+\mu 2^k\int_0^t{\zeta'(t') ||\Delta^h_k(\varepsilon^2\partial_x,\varepsilon\partial_y)(Q_{11},\varepsilon Q_{12})_\Theta)||^2_{L^2}dt'}\\
&\le ||e^{a|D_x|}\Delta^h_k(\varepsilon (u^\varepsilon_0-u_0), \varepsilon ^2(v^\varepsilon_0-v_0))||^2_{L^2}+\int_0^t{|(\Delta^h_k\mathcal S^1_\Theta, \Delta^h_k w^1_\Theta)_{L^2}|dt'}+\int_0^t{|(\Delta^h_k\mathcal S^2_\Theta, \Delta^h_k w^2_\Theta)_{L^2}|dt'}\\
\end{aligned}
\end{equation}
We deduce from the Propsition \ref{proplim1}, \ref{proplim2} and $M>1$ that
\begin{equation}
\label{dim2eqcon1}
\begin{aligned}
&\int_0^t{|(\Delta^h_k\mathcal S^1_\Theta, \Delta^h_k w^1_\Theta)_{L^2}|dt'}+\int_0^t{|(\Delta^h_k\mathcal S^2_\Theta, \Delta^h_k w^2_\Theta)_{L^2}|dt'}\\&\lesssim d_k^22^{-k}(||(\varepsilon w^1, \varepsilon^2 w^2)_\Theta||^2_{\tilde L^2_{t,\zeta'(t)}(B^1)}+M\varepsilon||(\varepsilon\partial_x,\partial_y)(\varepsilon w^1_\Theta, \varepsilon^2 w^2_\Theta)||_{\tilde L^2_t(B^\frac 12)})\\
&+d_k^22^{-k}(M^\frac 12 \varepsilon||\partial_y(\varepsilon w^1_\Theta, \varepsilon^2 w^2_\Theta)||_{\tilde L^2_t(B^\frac 12)}||(\varepsilon w^1_\Theta, \varepsilon^2 w^2_\Theta)||_{\tilde L^2_{t,\zeta'(t)}(B^ 1)}+M^\frac 32\varepsilon||\varepsilon^2 w^2_\Theta||_{\tilde L^2_{t,\zeta'(t)}(B^1)})\\
\end{aligned}
\end{equation}
so from Cauchy-Schwartz inequality, we have
\begin{equation}
 \begin{aligned}
 \label{dim2eqcon2}
 &||(\varepsilon w^1)_{\Theta}, \varepsilon^2 w^2_\Theta)||_{\tilde L_t^\infty(B^\frac 12)}+\mu ||(\varepsilon w^1)_{\Theta}, \varepsilon^2 w^2_\Theta)||_{\tilde L_{t,\zeta'(t)}^2(B^1)}+||(Q_{11}^{\varepsilon},\varepsilon Q^{\varepsilon}_{12})_\Theta||_{\tilde L^\infty_t(B^\frac 12)}\\
 &+||(\varepsilon^2\partial_x,\varepsilon\partial_y)(Q_{11}^{\varepsilon},\varepsilon Q^{\varepsilon}_{12})_\Theta||_{\tilde L^\infty_t(B^\frac 12)}+||\partial_y(\varepsilon w^1,\varepsilon^2w^2)_\Theta||_{\tilde L^2_t(B^\frac 12)}+||e^{\mathcal Rt'}\partial_x(\varepsilon u^{\varepsilon},\varepsilon^2v^{\varepsilon})_\psi||_{\tilde L^2_t(B^\frac 12)}\\
&+||\partial_y(Q^{\varepsilon}_{11},\varepsilon Q^{\varepsilon}_{12})_\Theta||_{\tilde L^2_t(B^\frac 12)}+||\varepsilon\partial_x(Q^{\varepsilon}_{11},\varepsilon Q^{\varepsilon}_{12})_\Theta||_{\tilde L^2_t(B^\frac 12)}+\varepsilon||\Delta_\varepsilon(Q^{\varepsilon}_{11},\varepsilon Q^{\varepsilon}_{12})_\psi||_{\tilde L^2_t(B^\frac 12)}\\
&\le C||e^{a|D_x|}(\varepsilon(u^\varepsilon_0- u_0),\varepsilon^2(v_0^\varepsilon-v_0))||_{B^\frac 12}+C||e^{a|D_x|}((Q_{11})_0,\varepsilon (Q_{12})_0)||_{B^\frac 12}\\
&+C||e^{a|D_x|}(\varepsilon^2\partial_x,\varepsilon\partial_y)((Q_{11})_0,\varepsilon (Q_{12})_0)||_{B^\frac 12}+C(||(\varepsilon w^1, \varepsilon^2 w^2)_\Theta||_{\tilde L^2_{t,\zeta'(t)}(B^1)}\\
&+M^\frac 12\varepsilon^\frac 12||(\varepsilon\partial_x,\partial_y)(\varepsilon w^1_\Theta, \varepsilon^2 w^2_\Theta)||^\frac 12_{\tilde L^2_t(B^\frac 12)}+M^\frac 34\varepsilon^\frac 12||\varepsilon^2w^2_\Theta||^\frac 12_{\tilde L^2_{t,\eta'(t)}(B^1)}\\
&+M^\frac 14 \varepsilon^\frac 12||\partial_y(\varepsilon w^1_\Theta, \varepsilon^2 w^2_\Theta)||^\frac 12_{\tilde L^2_t(B^\frac 12)}||(\varepsilon w^1_\Theta, \varepsilon^2 w^2_\Theta)||^\frac 12_{\tilde L^2_{t,\zeta'(t)}(B^ 1)}\\
&\le C||e^{a|D_x|}(\varepsilon(u^\varepsilon_0- u_0),\varepsilon^2(v_0^\varepsilon-v_0))||_{B^\frac 12}+C||e^{a|D_x|}((Q_{11})_0,\varepsilon (Q_{12})_0)||_{B^\frac 12}\\
&+C||e^{a|D_x|}(\varepsilon^2\partial_x,\varepsilon\partial_y)((Q_{11})_0,\varepsilon (Q_{12})_0)||_{B^\frac 12}+CM(\varepsilon+||(\varepsilon w^1)_{\Theta}, \varepsilon^2 w^2_\Theta)||_{\tilde L_{t,\zeta'(t)}^2(B^1)})\\
 \end{aligned}
 \end{equation}
 and \eqref{dim2eqcon1} is proved by choosing  $\mu\ge C^2M^2$.
\end{proof}


\bigskip{\bf Acknowledgments.} X.Li and A.Zarnescu have been partially supported by the Basque Government through the BERC 2022- 2025 program and by the Spanish State Research Agency through BCAM Severo Ochoa excellence accreditation SEV-2017-0718 and through project PID2020-114189RB-I00 funded by Agencia Estatal de Investigaci´on (PID2020- 114189RB-I00 / AEI / 10.13039/501100011033). A.Zarnescu has been  also partially supported by a grant of the Ministry
of Research, Innovation and Digitization, CNCS - UEFISCDI, project number PN-III-P4-PCE-2021-0921, within
PNCDI III. M.Paicu has been supported by Universit\'e de Bordeaux.


\end{document}